# Accuracy analysis and optimization of scale-independent third-order WENO-Z scheme with critical-point accuracy preservation

Yunchuan Wu[a†], Xiuzheng Cheng[a†], Yi Duan[b], Linsen Zhang[b], Qin Li[a‡], Pan Yan[a], Mengyu Wang[a]

a. School of Aerospace Engineering, Xiamen University, China, 361102
b. Science and Technology on Space Physics Laboratory, China Academy of Launch Vehicle Technology, Beijing 100076, China

**Abstract**. To address the order degradation at critical points in the WENO3-Z scheme, some improvements have been proposed (e.g., WENO-NP3, -F3, -NN3, and -PZ3), but these approaches generally fail to consider the occurrence of critical points at arbitrary positions within grid intervals, resulting in their inability to maintain third-order accuracy when a first-order critical point ($CP_1$) occurs. Also, most previous improved schemes suffer from a relatively large exponent $p$ of the ratio of global to local smoothness indicators, which adversely affects the numerical resolution. Concerning these limitations, introduced here is an accuracy-optimization lemma demonstrating that the accuracy of nonlinear weights can be enhanced providing that smoothness indicators satisfy specific conditions, thereby establishing a methodology for elevating the accuracy of nonlinear weights. Leveraging this lemma, a local smoothness indicator is constructed with error terms achieving second-order in smooth regions and fourth-order at $CP_1$, alongside a global smoothness indicator yielding fourth-order accuracy in smooth regions and fifth-order at $CP_1$, enabling the derivation of new nonlinear weights that meet accuracy requirements even when employing $p = 1$. Furthermore, a resolution-optimization lemma is proposed to analyze the relationship between parameters in local smoothness indicators and resolution. By integrating theoretical analysis with numerical practices, free parameters in non-normalized weights and local smoothness indicators are determined under the balance of numerical resolution and robustness, which leads to the development of WENO3-$Z_{ES4}$, a new WENO3-Z improvement that preserves the optimal order at $CP_1$ especially with $p = 1$. 1D and 2D validating tests show that compared with previous WENO3-Z improvements, the new scheme consistently achieves third-order in the case of $CP_1$ regardless of its position and exhibits good resolution (e.g., the resolving of the second peak and valley in the 1D Shu–Osher problem with only 240 grid points) as well as preferable robustness.

## 1 Introduction

In computational fluid dynamics, extensive use is made of WENO (weighted essentially non-oscillatory) schemes [1,2], especially the WENO-JS version [2]. However, as is well known, WENO-JS suffers from order degradation at critical points, thereby prompting improvements such as the mapping-function approach WENO-M [3] and WENO-Z [4–7]; the former recovers optimal order primarily via a mapping function, while the latter does so by formulating a global smoothness indicator ($\tau$) (sometimes also local smoothness indicators $\beta_k$) to derive nonlinear weights. Toward the fifth-order WENO5-JS and via the two approaches separately, Henrick et al. [3] and Borges et al. [4] proposed respective improvements, with notable enhancements in accuracy recovery and resolution.

---

† Yunchuan Wu and Xiuzheng Cheng are co-first authors
‡ Corresponding author, email: qin-li@vip.tom.com

While both WENO-M and WENO-Z merit attention, the focus herein is solely on further research into the latter. Given the wide-ranging application of third-order schemes in practical engineering, especially their robustness and efficiency via fewer sub-stencils, the present study is restricted to corresponding improvements using two sub-stencils. In previous work [4,8], it was indicated that concerning the difference between linear and nonlinear weights, namely $\omega_k - d_k$, WENO3-Z only attains $O(\Delta x)$ in smooth regions, rather than the $O(\Delta x^2)$ required by the sufficient condition to achieve optimal order. Moreover, it was found that WENO3-Z fails to achieve third-order accuracy at first-order critical points ($CP_1$, where $f' = 0, f''$ and $f''' \neq 0$) [8,9]. To meet the accuracy requirements at $CP_1$, serial improvements such as WENO-NP3 [9], -F3 [10], -NN3 [11], and -PZ3 [12] have been proposed. Regarding these improvements, a recent study [13] stated the following. 1) Although motivated by different considerations, the construction of $\tau$ in these schemes shares the same form $c(f_{i-1} - 2f_i + f_{i+1})^2$ (where $f$ denotes the variables) except for the different evaluations of $c$. 2) When $CP_1$ occurs on a half-node, WENO-NP3, -F3, and -NN3 fail to achieve third-order in the $L_\infty$-norm, and WENO-NN3 even fails to attain first-order accuracy. 3) When performing analysis and scheme construction, previous studies usually assumed that critical points occur on grid nodes ($x_j$), neglecting their occurrence within grid intervals of the stencils involved and so raising problems in the precision relationship and the ensuing constructions, e.g., the occurrence in WENO-NN3 [11]. 4) The aforementioned WENO3-Z improvements are scale-dependent, meaning that the computational results for the same problem would differ when using different length or variable scales, as shown by the examples in [13]. Note that when studying the order degradation at critical points in WENO5-Z, Don and Borges [6] analyzed how the exponent $p$ and the small parameter $\varepsilon$ affect numerical accuracy, thereby offering insights for improving WENO3-Z.

As we mentioned previously [14], a critical point can occur at any position within the dependent stencil, which may cause the accuracy relationship of the smoothness indicator (abbreviated as IS according to [2]) to change with that position. With this understanding, a systematic investigation was reported in [13], a byproduct of which was a modified version of the WENO-NN3 scheme; subsequent computation proved its recovery of optimal order at $CP_1$, thereby validating the analysis concerning arbitrary critical-point location (ACPL) mentioned above. The consequences of APCL [13,15] further reveal that when $CP_1$ occurs in the scale-independent WENO3-Z, an IS derived from three-point stencils cannot satisfy the conditions for optimal order recovery. Consequently, by means of stencil extension, we developed the WENO3-ZM and WENO3-Z$_{ES}$ [13] schemes with the conditions satisfied. WENO3-ZM expands a single stencil point and uses the mapping function from [14], while WENO3-Z$_{ES}$ incorporates two extended points. However, subsequent numerical experiments showed both schemes to be insufficiently robust. To address this, we proposed two improvements with enhanced resolution and robustness, namely WENO3-Z$_{ES2}$ and WENO3-Z$_{ES3}$ [15], by constructing new global ISs with lower-order accuracy and optimized local ISs. In the improvements [15], the nonlinear weights are $\alpha_k = d_k(1 + C_\alpha(\tau/(\beta_k + \varepsilon))^p)$, where $d_k$ are linear weights, $C_\alpha$ is a constant coefficient, $\varepsilon$ is a small constant to prevent division by zero, and $p = 2$. Notably, studies showed an inverse effect between $p$ and the numerical resolution [4,13], i.e., the larger the value of $p$, the less the resolution, which underscores the practical and theoretical importance of $p$ optimization. Likewise, such optimization should comply with accuracy requirements from ACPL, scale-independence, and the maintenance of numerical robustness.

Of note are the WENO variants WENO-P [16], WENO-AO [17], TENO [18], nonlinear BVD [19,20], and AI-WENO (affine-invariant WENO) [21]. In WENO-P [16], the ISs $\tilde{\beta}_k^{(3)}$ in the $L_1$-norm of three-point stencils were devised and a global $\zeta$ analogous to that in WENO5-Z was defined, with nonlinear weights formulated via $(\zeta/\tilde{\beta}_k^{(3)})^2$. WENO-AO [17] hybridizes high-order linear schemes with lower-order (at least third-order) linear schemes via combinations. TENO [18] achieves enhanced robustness by dynamically assembling incremental low-order stencils and leveraging linear combinations of low-order approximations. As a discontinuity-capturing method in the finite-volume framework, nonlinear BVD [19,20] minimizes boundary variation by adaptively selecting reconstruction functions, i.e., high-order polynomials in smooth regions and non-polynomial functions (e.g., THINC functions) near discontinuities. AI-WENO [21] introduces a descaler and average operation to preprocess primitive variables, enabling nonlinear weights to adapt to multiscale function variations and thereby improving robustness. However, these advancements primarily target schemes of fifth-order or higher, and investigations into third-order WENO schemes remain limited.

Building on the above, herein we propose a lemma regarding precision optimization and establish a corresponding approach to enhance the accuracy of normalized nonlinear weights. In the WENO-Z framework, we begin by formulating an IS via this lemma, followed by related parametric analysis and subsequent derivation of another lemma regarding how parameters in local ISs affect numerical resolution. Based on these outcomes and the results of numerical experiments, the parameters of the scheme (i.e., $C_\alpha$ and $C_\beta$) are determined, and an improved scheme is proposed with enhanced resolution, preserved robustness, and scale-independence. Specifically, Section 2 reviews the WENO-Z framework and classical third-order variants, introducing ACPL and relevant theoretical progress. Section 3 presents the lemmas and consequent scheme construction, with the parameters $C_\beta$ and $C_\alpha$ determined by analysis and simulations, ultimately yielding the new WENO3-Z$_{ES4}$. Section 4 reports comprehensive canonical tests using WENO3-Z$_{ES4}$, as well as comparison with WENO3-Z and two improvements regarding resolution and robustness. Finally, concluding remarks are made in Section 5.

## 2 Review of WENO3-Z, Relevant Analysis, and Improvements

Because the present study uses the WENO-Z framework, in this section we present a brief review of the scheme, especially WENO3-Z and its improvements, including analyses concerning APCL and theoretical outcomes.

### 2.1 WENO-Z Formulation

For convenience, we begin by revisiting WENO-JS [1,2]. We consider the 1D hyperbolic conservation law

$$u_t + f(u)_x = 0, \tag{1}$$

where the flux derivative satisfies $f(u)_x = \frac{\partial f(u)}{\partial x} > 0$. The grid points $x_j = j\Delta x$ are defined on a uniform grid with spacing $\Delta x$ and indexed by $j$. The conservative discretization of $f(u)_x$ at $x_j$ is formulated by the reconstruction $\hat{f}$ as

$$(f(u)_x)_j \approx (\hat{f}_{j+\frac{1}{2}} - \hat{f}_{j-\frac{1}{2}})/\Delta x. \tag{2}$$

In WENO-JS, we have

$$\hat{f}_{j+\frac{1}{2}} = \sum_{k=0}^{r-1} \omega_k q_k^r = \sum_{k=0}^{r-1} \omega_k \sum_{l=0}^{r-1} a_{kl}^r f(u_{(j-r+k+l+1)}), \tag{3}$$

where $q_k^r$ denotes the linear interpolation on the $k$-th candidate stencil, $a_{kl}^r$ are the interpolation coefficients, $\omega_k$ are the normalized nonlinear weights, and $r$ is the grid number of the candidate stencil. We have $\omega_k = \alpha_k / \sum_{l=0}^{r-1} \alpha_l$, where the non-normalized weight $\alpha_k$ from [2] is $d_k/(\varepsilon + \beta_k^{(r)})^2$ with $\varepsilon = 10^{-6} \sim 10^{-7}$ and $\beta_k^{(r)}$ as

$$\beta_k^{(r)} = \sum_{m=0}^{r-2} c_m^r \left( \sum_{l=0}^{r-1} b_{kml}^r f(u_{j-r+k+l+1}) \right)^2. \tag{4}$$

The coefficients $a_{kl}^r$, $c_m^r$, and $b_{kml}^r$ in $\beta_k^{(r)}$ are given in [13,14]. For the third-order case or $r = 2$, we have $\beta_0^{(2)} = (f_j - f_{j-1})^2$ and $\beta_1^{(2)} = (f_{j+1} - f_j)^2$. In WENO-Z [4–6], the nonlinear weights take the form

$$\alpha_k = d_k \left( 1 + C_\alpha \left( \tau/(\beta_k^{(r)} + \varepsilon) \right)^p \right), \tag{5}$$

where $\tau$ is the global IS, $p$ takes the value 1 or 2, $C_\alpha$ is a coherently free coefficient (typically set to 1), and $\varepsilon$ is extremely small ($10^{-40}$ or smaller).

## 2.2 APCL Analysis and Theoretical Outcomes

In this section, the analysis concerns WENO3-Z and similar schemes, the stencils of which span $\{x_{j-1}, x_j, x_{j+1}\}$. ACPL [13,14] refers to the fact that when a critical point occurs, it may exist at arbitrary positions within the stencil dependent on other than $x_j$, i.e., its coordinate is

$$x_c = x_j + \lambda \cdot \Delta x, -1 < \lambda < 1. \tag{6}$$

Expanding ISs at the critical point considering Eq. (6) would yield leading errors of both $\beta_k$ and $\tau$ that depend explicitly on $\lambda$ and probably differ from the conventional ones assuming a critical point at $x_j$ ($\lambda = 0$ only). For example, consider WENO3-Z with the global IS as $\tau_3 = |(f_{j+1} - f_{j-1})(f_{j+1} - 2f_j + f_{j-1})|$. At $CP_1$, the leading error of $\tau_3$ is $|-2\lambda f_{x_c}''^2 \Delta x^4|$, while that of $\beta_{0,1}^{(2)}$ is $\frac{1}{4}(2\lambda \pm 1)^2 f_{x_c}''^2 \Delta x^4$. Therefore, according to ACPL [13,14], $\tau_3$ achieves only fourth-order accuracy for $\lambda \neq 0$ while becoming fifth-order at $\lambda = 0$; similarly, $\beta_{0,1}^{(2)}$ retains fourth-order accuracy provided that $\lambda \neq \pm 1/2$, otherwise increasing to sixth-order. From [3,4], the sufficient condition for achieving the optimal $(2r - 1)$th-order of WENO is

$$\omega_k^\pm - d_k = O(\Delta x^r), \tag{7}$$

where the superscript "$\pm$" denotes positions at $x_{j\pm\frac{1}{2}}$, and $r$ has the same meaning as before.

Subsequent derivations confirm that the sufficient condition for WENO3-Z achieving third-order accuracy should be $(\tau/\beta_k + \varepsilon))^p = O(\Delta x^2)$, and it is known from the above analysis that WENO3-Z fails to preserve third-order accuracy in the presence of critical points. To handle this critical-point degradation, new IS and nonlinear weighting strategies must be designed. As shown later, Eq. (9) yields nonlinear weights to likely achieve optimal order at critical points, but its scale-dependent formulation introduces non-unique solutions when applied across differing computational scales. To construct a scale-independent scheme, the nonlinear weights should instead be computed using $d_k(1 + (\tau/\beta_k + \varepsilon))^p)$. Notably, because the dimension of $\beta_k$ is $[f]^2$, $\tau$ must maintain the same dimension to ensure scale-independence. To construct an appropriate global IS, we proved the following lemma in [13].

**Lemma 2.1:** Consider the three-point stencil $\{x_{j-1}, x_j, x_{j+1}\}$. When constructing a global smoothness indicator $\tau(f)$ as a quadratic function of $f$ in the form

$$\tau(f) = (f_{j-1}, f_j, f_{j+1})[a_{i_1,i_2}](f_{j-1}, f_j, f_{j+1})^T, \tag{8}$$

where $[a_{i_1,i_2}]$ is a $3 \times 3$ matrix, the following conclusions hold:
(1) Without critical points, the only formulation allowing $\tau(f)$ to expand to $O(\Delta x^4)$ is $\tau(f) = c(f_{j-1} - 2f_j + f_{j+1})^2$.

(2) With $CP_1$ ($f' = 0, f''$ and $f''' \neq 0$) occurring at $x_c = x_j + \lambda \cdot \Delta x$ ($-1 < \lambda < 1$), there exist no non-trivial solutions $a_{i1,i2}$ such that $\tau(f)$ expands to $O(\Delta x^5)$ at the critical point.

From Lemma 2.1(1), although $\tau(f)$ in typical WENO-Z improvements is constructed from different motivations, each fourth-order $\tau(f)$ must have the same form with only distinct evaluations of $c$ (see Table 1 in Section 3 for specific schemes and coefficients). Lemma 2.1(2) indicates that when $CP_1$ occurs at an arbitrary position within the interval ($-1 < \lambda < 1$), the stencil $\{x_{j-1}, x_j, x_{j+1}\}$ cannot produce a scale-independent scheme with third-order accuracy unless stencil extension is performed. Based on this understanding, we developed the improved WENO3-Z schemes in [13,15], as detailed in Section 3. To clarify the influence of coefficients on resolution, we further proved the following lemma [13].

**Lemma 2.2:** Consider two candidates $S_C$ and $S_D$ corresponding to the same stencil, where "$C$" and "$D$" indicate that the variable is distributed more smoothly on $S_C$ than on $S_D$ or $0 < \beta_{k,C} < \beta_{k,D}$.
(1) For nonlinear weights $\alpha_k = d_k(1 + (\tau/\beta_k)^{p_i})$ with $p_i > 0$ and $\tau > 0$, if $p_1 > p_2 > 0$ and $0.278 < \tau/\beta_{k,D} < \tau/\beta_{k,C}$, then $[\omega_{k,D}/\omega_{k,C}]_{p_1} < [\omega_{k,D}/\omega_{k,C}]_{p_2}$.

(2) For nonlinear weights $\alpha_k = d_k(1 + C_{\alpha i}(\tau/\beta_k)^p)$ with $C_{\alpha i} > 0$, $\tau > 0$ and $p > 0$, if $C_{\alpha 1} > C_{\alpha 2}$, then $[\omega_{k,D}/\omega_{k,C}]_{C_{\alpha 1}} < [\omega_{k,D}/\omega_{k,C}]_{C_{\alpha 2}}$.

Borges et al. [4] noted that if the conditions for achieving optimal order are satisfied, then a larger relative importance of the discontinuous stencil increases the numerical resolution. From Lemma 2.2(1), a larger (resp. smaller) $p$ decreases (resp. increases) the resolution, so to balance resolution and robustness, further parameter tunings are required. Lemma 2.2(2) shows that a smaller (resp. larger) $C_\alpha$ increases (resp. decreases) the resolution.

In summary, particularly from Lemma 2.1(2), achieving optimal order recovery at critical points requires either extending the stencil to construct scale-independent schemes or retaining the existing stencil with scale-dependent schemes available. In Section 2.3, we introduce representative WENO3-Z improvements.

## 2.3 Typical WENO3-Z Improvements to Recover Optimal Order at Critical Points

The usual WENO3-Z improvements assume that critical points occur only at $x_j$, with unnormalized weights having the form

$$\alpha_k = d_k(1 + \tau^{p_1}/(\beta_k + \varepsilon)^{p_2}), \tag{9}$$

where $p_1 \neq p_2$ generally. The global IS typically has one of the following forms:

$$\tau = \begin{cases} c_{\tau_1}|(f_{j+1} - f_{j-1})(f_{j+1} - 2f_j + f_{j-1})| & \text{or} \\ c_{\tau_2}(f_{j+1} - 2f_j + f_{j-1})^2, & \end{cases} \tag{10.a}$$
$$\tag{10.b}$$

where the values of $c_{\tau_{1,2}}$ differ for different schemes. By choosing respective $c_{\tau_{1,2}}$ and $p_{1,2}$, various improvements can be derived (see the specific choices in Table 1). Notably, WENO-NP3, -F3, and -NN3 use the $\tau$ of Eq. (2.10.b), and they all assume that critical points occur at $x_j$ when performing accuracy analysis at $CP_1$. ACPL reveals that when a critical point appears at a half-node, the optimal order cannot be achieved by the above improvements. Specifically, WENO-NN3 may even fail to attain first-order accuracy in cases of critical points. To address this, [13] provided new $p_1$ and $p_2$ for WENO-NN3 based on ACPL, via which it was validated that optimal order is recovered.

Table 1: Parameters in Eqs. (9) and (10) for typical improvements to WENO3-Z.

| Scheme | $\alpha_k$ | | $\tau$ | |
|---|---|---|---|---|
| | $p_1$ | $p_2$ | $c_{\tau_1}$ | $c_{\tau_2}$ |
| WENO-NP3 [9] | 3/2 | 1 | – | 10/12 |
| WENO-F3 [10] | 3/2 | 1 | – | 2/12 |
| WENO-NN3 [11] | 1 | ≤ 3/4 | – | 10/12 |
| WENO-PZ3 [12] | 1 | ≤ 1/2 | 1 | – |

As noted in [13], all the improvements in Table 1 are scale-dependent (with $p_1 \neq p_2$), leading to the anomalous situation of the same scheme giving different solutions to the same problem when different scales are used. Therefore, it is necessary and important to develop scale-independent improvements. Based on the analysis in Section 2.2, it is essential to extend the stencil to achieve third-order accuracy when critical points occur. In [13], a high-order global IS was constructed as

$\tau_{CP_1} = c \times (-f_{j+2} + 3f_{j+1} + 21f_j - 23f_{j-1}) \times (f_{j+2} - 3f_{j+1} + 3f_j - f_{j-1})$, and $\beta_1^{(3)}$ on the downwind stencil of WENO5-JS was used to replace $\beta_1^{(2)}$. Combined with the mapping-function method, this yielded the improved WENO3-ZM, which preserves optimal order at $CP_1$. Also proposed was a global IS with even higher order, i.e., $\tau_{CP_2} = c(f_{j+2} - 4f_{j+1} + 6f_j - 4f_{j-1} + 4f_{j-2})^2$ integrated with the local ISs $\beta_0^{(3)}$ and $\beta_1^{(3)}$. This resulted in the WENO3-Z$_{ES}$ scheme, which maintains third-order accuracy even at a second-order critical point. However, although WENO3-ZM and WENO3-Z$_{ES}$ achieve third-order accuracy at $CP_1$, their robustness remains limited. To address this, in [15] we developed two more stable schemes, i.e., WENO3-$Z_{ES2}$ and WENO3-$Z_{ES3}$. These schemes use a four-point stencil to construct a new global IS and two types of local ISs. The two schemes differ in their local ISs on the first candidate stencil; i.e., in WENO3-$Z_{ES2}$, we have $\beta_0^{(2)*} = \beta_0^{(2)} + C_{\beta 0}(\delta_j^{(2)_2})^2$, where $C_{\beta 0}$ is determined via discontinuity detection, and $\delta_j^{(2)_2} = (f_{j+1} - 2f_j + f_{j-1})$. In WENO3-$Z_{ES3}$, we have $\beta_0^{(2)*} = \beta_0^{(2)} + C_{\beta 0}(\delta_{j-1}^{(2)_2})^2$, where $C_{\beta 0} = 0.5$ and $\delta_{j-1}^{(2)_2} = (f_j - 2f_{j-1} + f_{j-2})$. Both schemes use the same $\beta_1^{(2)*} = \beta_1^{(2)} + C_{\beta 1}(\delta_{j+1}^{(2)_2})^2$, where $\delta_{j+1}^{(2)_2} = (f_{j+2} - 2f_{j+1} + f_j)$ and $C_{\beta 1} = 0.15$. For discussion purposes, $\delta_j^{(n)m}$ denotes below the order-$n$ derivative difference with order-$m$ accuracy at point $x_j$. Also, a $\tau$ is constructed for two improvements on $\{x_{j-1}, x_j, x_{j+1}, x_{j+2}\}$ attaining $O(\Delta x^5)$ at $CP_1$, i.e., $\tau = |(2f_{j+1} - 3f_j + f_{j-1})(2f_{j+2} - 3f_{j+1} + 3f_j - f_{j-1})|$. The scale-independent nonlinear weights are formulated as

$$\alpha_k = d_k(1 + C_\alpha(\tau/(\beta_k + \varepsilon))^p), \tag{11}$$

where $p = 2$ in terms of the sufficient condition Eq. (7). From Lemma 2.2, if a $p$ less than the normally-used 2 can make the scheme achieve third-order even at $CP_1$, then the numerical resolution would be enhanced. To achieve this, we propose an accuracy-optimization lemma based on which we then develop WENO3-$Z_{ES4}$, a scale-independent improved scheme with balanced resolution and robustness; see Section 3 for details.

## 3 Improving WENO3-Z via Accuracy Optimization Analysis

In Section 2, we briefly reviewed WENO3-Z and its typical improvements. In this section, we investigate accuracy optimization under the WENO3-Z framework, beginning by introducing two lemmas related to scheme improvement.

### 3.1 Lemmas for Optimization of Accuracy and Resolution

The analysis in Section 2 shows that in our previous improvements, i.e., WENO3-ZM, -$Z_{ES2}$, and -$Z_{ES3}$, $\tau/\beta_k$ is mostly of $O(\Delta x)$, which requires the exponent $p$ to be greater than 1 and may consequently decrease numerical resolution. Naturally, the question is whether it is possible to construct a scale-independent improvement by $p = 1$ with optimal order at $CP_1$ achieved and favorable robustness maintained. A straightforward approach is to construct a high-order global IS to achieve $\tau/\beta_k = O(\Delta x^2)$ when $CP_1$ occurs. However, we have shown that this would make $\tau/\beta_k$ have a rather high order in the absence of critical points, which empirically compromises the numerical robustness. Therefore, we aim to design a scheme that upgrades the order of the nonlinear weights, enabling them to attain the second-order required by Eq. (7) when $CP_1$ occurs, through which $\tau$ would have an order that is not too high. To achieve this goal, in-depth analysis is imperative, the outcome of which leads to the following lemma.

**Lemma 3.1:** For normalized nonlinear weights $\omega_k = \alpha_k/\sum_l \alpha_l$ where $\alpha_k = d_k(1 + C_{\alpha k}(\tau/\beta_k)^p)$ and $\tau = O(\Delta x^m)$, suppose that the accuracy relationship of $\beta_k$ satisfies

$$\beta_k = \sum_{l'=n_1}^{n_2-1} a_{l'}\Delta x^{l'} + b_k \Delta x^{n_2} + O(\Delta x^{n_2+1}), \tag{12}$$

where $b_k$ varies across different $k$, $n_2 > n_1 \geq 1$, and typically $m > n_1$. If $C_{\alpha k}$ take the same value for different $\alpha_k$, then $\omega_k = d_k(1 + O(\Delta x^{p(m-n_1)}) \times O(\Delta x^{n_2-n_1}))$.

Lemma 3.1 is proved in the Appendix. The usual understanding is that $\omega_k$ is of $d_k(1 + O(\Delta x^{p(m-n_1)}))$, but Lemma 3.1 indicates that the order of $\omega_k$ is upgraded by $(n_2 - n_1)$ compared to the original $p(m - n_1)$. By this lemma, we can construct improved schemes with $p = 1$. Specifically, for schemes targeting the optimal $(2r - 1)$th-order accuracy, we select appropriate local ISs such that their precision relationships at critical points satisfy the conditions by Eq. (12), i.e., coefficients of $\Delta x^{l'}$ with $l' \leq (n_2 - 1 - n_1)$ are identical, while those of terms with higher order differ. Because the accuracy of $\omega_k$ is upgraded by $n_2 - n_1$, the requirement on the magnitude of $\tau$ is relaxed from the conventional $O(\Delta x^{r+n_1})$ to $O(\Delta x^{r+n_1-(n_2-n_1)})$. Note that the $C_{\alpha k}$ of $\alpha_k$ must be the same, differing from the canonical WENO-Z in which $C_{\alpha k}$ may vary. For clarity, we shorten $C_{\alpha k}$ to $C_\alpha$.

In Section 2.3, we noted that in WENO3-$Z_{ES2}$ and -$Z_{ES3}$, the $\delta_j^{(2)_2}$ in $\beta_k^{(2)*}$ has the coefficient $C_\beta$. Our previous numerical experiments showed that larger $C_\beta$ enhances numerical resolution while compromising robustness. It is well-established that a higher relative importance of the discontinuous stencil is apt to improve resolution, and based on this we propose the following lemma.

**Lemma 3.2:** Consider two candidates $S_C$ and $S_D$ corresponding to the same stencil, where "$C$" and "$D$" denote that the variable is smoother on $S_C$ than on $S_D$ or $\beta_{k,D} > \beta_{k,C}$. Suppose $\beta_k = (\delta_j^{(1)m})^2 + C_{\beta,i}(\delta_j^{(2)n})^2$, where $\delta_j^{(n)m}$ represent the $n$-th order derivative with $m$-th order accuracy at point $x_j$, and the subscript $i$ indicates the different evaluation of $C_{\beta,i}$. Likewise, $\delta_{j,D}^{(1)m} > \delta_{j,C}^{(1)m}$ and $\delta_{j,D}^{(2)n} > \delta_{j,C}^{(2)n}$. For $\alpha_k = d_k(1 + \tau/\beta_k)$ with $\tau > 0$ and $\delta_{j,D}^{(1)m^2}/\delta_{j,C}^{(1)m^2} > \delta_{j,D}^{(2)n^2}/\delta_{j,C}^{(2)n^2}$, then $[\omega_{k,D}/\omega_{k,C}]_{C_{\beta,1}} > [\omega_{k,D}/\omega_{k,C}]_{C_{\beta,2}}$ providing that $C_{\beta,1} > C_{\beta,2} > 0$.

Lemma 3.2 is proved in the Appendix. From this lemma, it is evident that a larger $C_\beta$ generally increases the relative importance of the discontinuous stencil, thus numerically enhancing resolution; conversely, reducing $C_\beta$ favors robustness. Recalling Lemma 2.2, theoretical findings suggest that the adjustment potential of $C_\beta$ and $C_\alpha$ would enhance the resolution while maintaining robustness.

## 3.2 Improving WENO3-Z via Accuracy and Resolution Optimization

In this subsection, we construct appropriate ISs so as to apply the lemmas proposed above, further determine the coefficients in the scheme integrated with theoretical analysis and numerical practices, and finally present the new WENO3-Z improvement named WENO3-$Z_{ES4}$.

### 3.2.1 Construction of ISs Considering Accuracy-optimization Lemma

For third-order schemes, the emergence of critical points induces an order increase for the leading error of $\beta_k$. To preserve optimal order at critical points, the leading error of $\tau$ usually has high order, which is apt to satisfy the sufficiency conditions by Eq. (7) in the absence of critical points. Therefore, current analysis prioritizes accuracy relationships in the case of $CP_1$. To apply Lemma 3.1, the local ISs must be formulated to satisfy Eq. (12). Our past computations revealed empirically that $\tau$ having higher order may compromise numerical robustness. To ensure as much robustness as possible, efforts should be made to decrease the order of $\tau$, which consequently reduces those of $\beta_k$. In [13], order increase of $\beta_k$ when $CP_1$ occurs at a half-node was prevented by expanding the stencil, i.e., using $\beta_2^{(3)} = \frac{1}{4}(3f_j - 4f_{j+1} + f_{j+2})^2 + \frac{13}{12}(f_j - 2f_{j+1} + f_{j+2})^2$ other than $\beta_1^{(2)}$. In $\beta_2^{(3)}$, the first term represents the squared first-derivative approximation at $x_j$ and the second term corresponds to the squared second-derivative approximation. Our subsequent investigations [13] demonstrated numerically the improved robustness when replacing the original three-point approximation of the first-order derivative with a two-point one. Building on this, [13] derived the modified $\beta_1^{(2)*}$ as well as $\beta_0^{(2)*}$, i.e.,

$$\beta_0^{(2)*} = (f_j - f_{j-1})^2 + c_{\beta_0}(f_{j-2} - 2f_{j-1} + f_j)^2, \tag{13}$$

$$\beta_1^{(2)*} = (f_{j+1} - f_j)^2 + c_{\beta_1}(f_j - 2f_{j+1} + f_{j+2})^2. \tag{14}$$

Expanding these at $CP_1$ yields

$$\beta_0^{(2)*} = (f'''^2 \lambda^2 + f'''^2 \lambda + f'''^2/4 + f'''^2 c_{\beta_0})\Delta x^4 + O(\Delta x^5), \tag{15}$$

$$\beta_1^{(2)*} = (f'''^2 \lambda^2 - f'''^2 \lambda + f'''^2/4 + f'''^2 c_{\beta_1})\Delta x^4 + O(\Delta x^5). \tag{16}$$

It is evident that regardless of evaluations of $c_{\beta_0}$ and $c_{\beta_1}$, the formulations fail to fulfill Eq. (12). In terms of Lemma 3.1, the following local ISs similar to $\beta_k^{(3)}$ are adopted:

$$\beta_0^* = \frac{1}{4}(3f_j - 4f_{j-1} + f_{j-2})^2 + C_{\beta_0}(f_{j-2} - 2f_{j-1} + f_j)^2, \tag{17}$$

$$\beta_1^* = \frac{1}{4}(3f_j - 4f_{j+1} + f_{j+2})^2 + C_{\beta_1}(f_j - 2f_{j+1} + f_{j+2})^2, \tag{18}$$

where $C_{\beta_k} > 0$ and the stencil spans $\{x_{j-2}, x_{j-1}, x_j, x_{j+1}, x_{j+2}\}$. In the absence of critical points, $\beta_0^*$ and $\beta_1^*$ exhibit magnitudes of $O(\Delta x^2)$. When $CP_1$ occurs, their accuracy relations become

$$\beta_0^* = (\lambda^2 + C_{\beta 0})f'''^2 dx^4 + b_1 dx^5 + O(\Delta x^6), \tag{19}$$

$$\beta_1^* = (\lambda^2 + C_{\beta 1})f'''^2 dx^4 + b_2 dx^5 + O(\Delta x^6), \tag{20}$$

where $-1 < \lambda < 1$. At $CP_1$, the local ISs attain the magnitude of $O(\Delta x^4)$. To have $\beta_k^*$ satisfy the accuracy relationship $\beta_k = \sum_{l'=n_1}^{n_2-1} a_{l'} \Delta x^{l'} + b_k \Delta x^{n_2} + O(\Delta x^{n_2+1})$ (see Lemma 3.1), $C_{\beta_k}$ must follow $C_{\beta 0} = C_{\beta 1}$. For clarity, we shorten them to $C_\beta$ with the value determined in Section 3.2.2. Specifically, for formulations by Eqs. (19) and (20), we have $n_1 = 4$ and $n_2 = 5$ for those in

Eq. (12). By Lemma 3.1, the accuracy upgrade of $\omega_k$ is $(n_2 - n_1) = 1$. Notably in Eq. (12), we have $b_k = (-4\lambda^3 + \frac{2}{3}\lambda - 2C_\beta\lambda \pm 2C_\beta)$. When $-1 < \lambda < 1$, $b_k$ vanishes at specific points, which would increase $n_2$ and further elevate the error order of $\omega_k$. However, this occurrence does not affect the scheme's formal accuracy, and therefore the use of Eqs. (19) and (20) is retained.

Via the preceding analysis, the local ISs have been defined. Evidently, $O(\Delta x^2)$ is maintained in the absence of critical points, reducing to $O(\Delta x^4)$ when $CP_1$ emerges. For third-order schemes ($r = 2$), according to Lemma 3.1, to recover the optimal order upon $CP_1$, the corresponding $\tau$ for $\beta_k = \beta_k^*$ satisfying Eq. (12) should be of $O(\Delta x^{r+n_1-(n_2-n_1)}) = O(\Delta x^5)$. Considering the requirement of scale-independence, $\tau$ must also be a quadratic function of $f$. In [14], it was shown that $\tau$ essentially represents the product of two difference approximations of derivatives at $x_j$. Therefore, to obtain a $\tau$ of $O(\Delta x^5)$ at $CP_1$, at least a difference approximation of a third-order derivative is required. This yields a global IS candidate as $\tau = |\delta_j^{(1)} \delta_j^{(3)}|$. Specifically, [15] provided a global IS with favorable mathematical properties, i.e.,

$$\tau_4 = |(f_{j+2} - 3f_{j+1} + 3f_j - f_{j-1})(2f_{j+1} - 3f_j + f_{j-1})|. \tag{21}$$

The error of $\tau_4$ is $|f'f^{(3)}\Delta x^4 + O(\Delta x^5)|$ in the absence of critical points and becomes $|(3/2 - \lambda)f''_{xc}f^{(3)}_{xc}\Delta x^5 + O(\Delta x^6)|$ when $CP_1$ emerges. Under the situation of no critical points, because $\tau_4/\beta_k = O(\Delta x^2)$, the optimal order is certainly achieved.

### 3.2.2 Determination of $C_\beta$ in $\beta_k^*$ and $C_\alpha$ in $\alpha_k$

In Section 3.2.1, the fundamental structure of the scheme was established via the construction of ISs. As is well known, WENO-Z-type schemes contain free parameters such as $C_\beta$ and $C_\alpha$ that must be determined. As mentioned before, Lemma 3.2 shows that larger $C_\beta$ benefits resolution enhancement, and Lemma 2.2 shows that larger $C_\alpha$ degrades the resolution whereas smaller $C_\alpha$ improves it. Considering these lemmas, we determine $C_\beta$ and $C_\alpha$ via numerical experiments. Specifically, we use the Shu–Osher problem, which is a canonical indicator for numerical resolution (as detailed in Section 4). Through comparative numerical studies, the optimal values of $C_\beta$ and $C_\alpha$ are determined. Fig. 1 presents the density profiles under varying parameters, where "EXACT" denotes results obtained by using WENO5-JS on 10,000 grid points. The outcomes demonstrate enhanced resolution at both wave peaks and valleys with increasing $C_\beta$, while the increase of $C_\alpha$ degrades the resolution at these locations, which is consistent with the previous theoretical analysis.

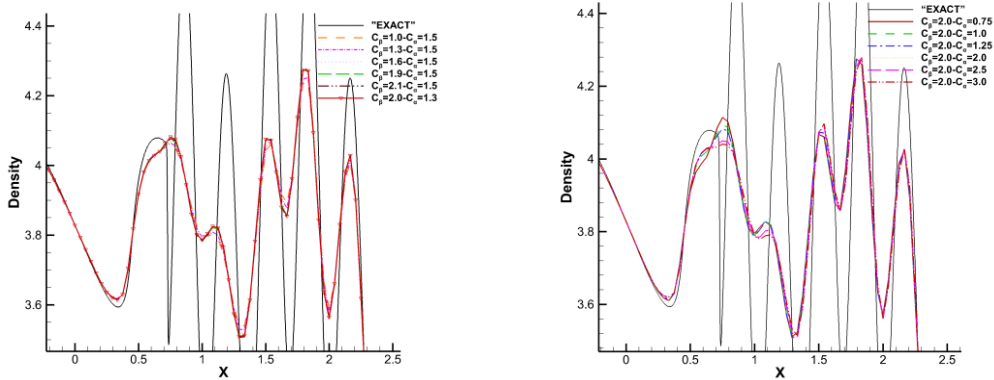

(a) $C_\alpha = 1.5$, $C_\beta = 1.0 \sim 2.1$ and $C_\alpha = 1.3$, (b) $C_\beta = 2.0$, $C_\alpha = 0.75 \sim 3.0$
$C_\beta = 2.0$

Figure 1: Density distributions for Shu–Osher problem by using different values of $C_\beta$ and $C_\alpha$ in WENO-$Z_{ES4}$ on 240 grid points at $t = 1.8$.

Note the following points. 1) Given the contradiction between numerical resolution and robustness, the determination of $C_\beta$ and $C_\alpha$ should balance the two attribute performances. 2) For the Shu–Osher problem, typical schemes barely resolve the second peak, making its resolution a pivotal criterion for parameter determination. As a tradeoff between resolution and robustness, we use $C_\beta = 2.0$ and $C_\alpha = 1.3$, and Fig. 1(a) shows the results obtained by using these values.

In summary, based on Lemma 3.1, we acquire a new scale-independent, third-order WENO3-Z improvement, herein called WENO3-$Z_{ES4}$. For completeness, its formulation is summarized as follows.

1) Local ISs $\beta_k^*$ are computed using Eqs. (17) and (18) with $C_{\beta 0} = C_{\beta 1} = 2.0$.

2) The global IS $\tau$ is derived via Eq. (21), then $\alpha_k$ are calculated by Eq. (5) using $C_\alpha = 1.3$ and $p = 1$.

3) $\alpha_k$ are normalized to yield $\omega_k$, which gives the reconstructed flux $\hat{f}_{j+\frac{1}{2}}$ by Eq. (3).

## 4. Validation Tests

### 4.1 Case Descriptions

For validation testing, we select the 1D scalar equation, the 1D and 2D Euler equations, and the 2D Navier–Stokes (N-S) equations, each chosen carefully as an indicator of numerical resolution and robustness. As well as the proposed WENO3-$Z_{ES4}$, the computations also use WENO3-Z ($p = 2$) [4–6], WENO3-F3 [10], and WENO3-$Z_{ES3}$ [15] for comparison. A concise overview of the case configurations is provided below.

1) 1D scalar equation

The equation is $\frac{\partial u}{\partial t} + \frac{\partial u}{\partial x} = 0$ with the initial condition of $u(x, 0)$. The region comprises $x \in [-1,1]$ with periodic boundary conditions. The following specific condition is selected:

$$u(x,0) = \sin\left(\pi(x - x_c) - \frac{\sin(\pi(x - x_c))}{\pi}\right), \tag{22}$$

where $x_c = 0.5966831869112089637212$. Two $CP_1$s occur initially at $x = 0$ and $x = -2 + 2x_c$. Temporal discretization uses the fourth-order Runge–Kutta (RK4) scheme with the time step $\Delta t = CFL \cdot \Delta x$. A series of refined grids is chosen with numbers $\{10, 20, 40, 80, \dots\}$, ensuring that $x = 0$ coincides with a grid node initially. A Courant–Friedrichs–Lewy (CFL) number of 0.25 is used to advance the computations to $t = 2$. Note that $CP_1$ located initially at $x = 0$ migrates to a half-node every four iterations. To the best of our knowledge, all existing WENO3-Z improvements without considering ACPL fail to achieve third-order accuracy in the $L_\infty$-norm.

2) 1D Euler equations

For the 1D Euler equations, three cases are selected: (i) the Shu–Osher problem, (ii) a blast wave, and (iii) a strong shockwave. For each case, the temporal scheme uses the third-order total variation diminishing (TVD) Runge–Kutta scheme (TVD-RK3), while the flux splitting uses the Steger–Warming scheme. In the conventional manner, characteristic variables are used to mitigate numerical oscillations.

(a) Shu–Osher problem

The initial conditions are $(\rho, u, p) = \{(3.857143, 2.62936, 10.3333), -5 \leq x < -4; (1 + 0.2\sin(5x), 0, 1), -4 < x \leq 5\}$. The computations are advanced to $t = 1.8$ on a uniform 240 grid points with $\Delta t = 0.003$. The results of WENO5-JS on 10,000 grid points are chosen as the "EXACT" solution.

(b) Blast wave

The initial conditions are $(\rho, u, p) = \{(1,0,1000), 0 \leq x < 0.1; (1,0,0.01), 0.1 \leq x \leq 0.9; (1,0,100), 0.9 < x \leq 1\}$ with solid-wall boundary conditions imposed on both ends. The computations are advanced to $t = 0.038$ on 600 grid points with $\Delta t = 1.0 \times 10^{-5}$. The results of WENO5-JS on 15,000 grid points serve as the "EXACT" solution.

(c) Strong shockwave

The initial conditions are $(\rho, u, p) = \{(1,0,0.1PR), -5 \leq x < 0; (1,0,0.1), 0 \leq x \leq 5\}$ where $PR = 10^6$. The computations are advanced to $t = 0.01$ on 200 grid points with $\Delta t = 1.0 \times 10^{-5}$.

3) 2D Euler equations

Six representative cases are chosen for the 2D Euler equations, including the 2D Riemann problem, double Mach reflection, shock–bubble interaction, and Mach-2000 jet flow. The temporal scheme and flux splitting method are the same as those for the 1D Euler equations.

(a) 2D Riemann problem

The domain spans $[0,1] \times [0,1]$ with the following initial conditions:

$$(\rho, u, v, p) = \begin{cases} (1.5, 0, 0, 1.5), 0.8 \leq x \leq 1, 0.8 \leq y \leq 1, \\ (0.5323, 1.206, 0, 0.3), 0 \leq x < 0.8, 0.8 \leq y \leq 1, \\ (0.138, 1.206, 1.206, 0.029), 0.8 \leq x < 0.8, 0 \leq y < 0.8, \\ (0.5323, 1.206, 0, 0.3), 0.8 \leq x < 0.8, 0 \leq y < 0.8. \end{cases}$$

The computations advance to $t = 0.8$ with $\Delta t = 0.0001$ on a $960 \times 960$ grid with a specific heat ratio of $\gamma = 1.4$.

(b) Double Mach reflection

This problem describes a Mach-10 shock impinging on a wall with a 60° incidence angle. The domain spans $[0,4] \times [0,1]$ with a 1920×480 grid. The initial conditions are $(\rho, u, v, p) = \{(8, 7.154, -4.125, 116.5), x \leq 1/6 + y/\sqrt{3}; (1.4, 0, 0, 1), x \geq 1/6 + y/\sqrt{3}\}$ with $\gamma = 1.4$. The computations are advanced to $t = 0.2$ with $\Delta t = 0.0001$.

(c) Shock–bubble interaction [22]

Also termed the Richtmyer–Meshkov instability, the shock–bubble interaction describes a

Mach-1.23 shockwave impinging on a density-stratified bubble where the internal density exceeds the outside density. The domain spans $[0,5] \times [0,5]$ with a shock initialized at $x = 4.5$ and a bubble centered at $(3,2.5)$ with radius $R = 1$. The initial conditions are as follows: $(\rho, u, v, p) = (0.2825, 0, 0.4721, 0.166)$ inside the bubble, $(\rho, u, v, p) = (0.2825, 0, 0.4721, 1)$ outside the bubble and pre-shock, and $(\rho, u, v, p) = (0, 0, 0.7546, 1.394)$ post-shock. The computations advance to $t = 10$ with $\Delta t = 0.001$ on a $1000 \times 1000$ grid.

(d) Mach-3 tunnel with a step [23]

This problem investigates a Mach-3 inviscid flow over a front step in a 2D tunnel, with the domain spanning $[0,5] \times [0,5]$. The step is located at $(0.6, 0)$ with a vertical height of 0.2. The initial flow conditions are $(\rho, u, v, p) = (1, 0, 1, 1)$ with $\gamma = 1.4$. The computations run until $t = 1.2$ at $\Delta t = 0.02$ on a $1200 \times 400$ grid.

(e) Inviscid sharp-double-cone flow at $Ma_\infty = 9.59$

This problem investigates high speed flow over a sharp double-cone configuration with semi-apex angles of 22° and 55°. The first geometric inflection point is positioned at $x = 3.63$, followed by a secondary corner at $x = 6$. A forward protrusion of length $L = 0.25$ extends from the cone vertex. The grid is a $204 \times 48$ one.

(f) Mach-2000 jet-flow problem [24]

As a canonical problem to test numerical robustness in extreme high speed flows, the domain spans $[0,1] \times [-0.25, 0.25]$ with the initial conditions $(\rho, u, v, p) = (0.5, 0, 0, 0.4127)$. An 800×400 grid is used, with outflow boundary conditions imposed on the right and lateral boundaries. For the left boundary, the conditions are as follows: $(\rho, u, v, p) = (5, 800, 0, 0.4127)$ for $|y| < 0.05$, otherwise $(\rho, u, v, p) = (0.5, 0, 0, 0.4127)$. Numerical integration proceeds with a time step of $\Delta t = 0.3 \times 10^{-7}$, advancing the solution to $t = 10^{-3}$.

4) Navier–Stokes equations

Three cases are selected: (i) flat-plate shock–boundary-layer interaction, (ii) the viscous shock-tube problem, and (iii) viscous sharp-double-cone flow at $Ma_\infty = 9.59$.

(a) Flat-plate shock–boundary-layer interaction at $Ma_\infty = 2$

The inflow conditions are $T_\infty = 293K$, $M_{a_\infty} = 2.0$, and $Re_\infty = 2.96 \times 10^5$ with an incident shock of angle 32.585°. The domain spans $[0, 2.02] \times [0, 1.3]$ discretized using a $103 \times 122$ grid. The boundary condition at the exit is supersonic outflow, and that on the surface is an adiabatic wall. The temporal discretization uses TVD-RK3.

(b) Viscous shock-tube problem [25]

This case provides complex flow structures of vortices and bifurcated shocks engendered by shock–boundary-layer interaction in a shock tube. The initial conditions are $(\rho, u, v, p) = \{(120, 0, 0, 120/\gamma), 0 \le x < 1/2; (1.2, 0, 0, 1.2/\gamma), 1/2 \le x < 1\}$ with $\gamma = 1.4$ The domain spans $[0,1] \times [0,1]$, and a $300 \times 300$ grid is used with no-slip and adiabatic wall boundaries on all sides. Key parameters include $Pr = 0.73$ and $Re = 200$. The computations advance to $t = 1$ with $\Delta t = 0.000416667$ using TVD-RK3 as the temporal scheme.

(c) Viscous sharp-double-cone flow at $Ma_\infty = 9.59$

This case serves as a critical test to assess the capabilities of computational fluid dynamics in predicting shock–boundary-layer interaction, and it has the same geometry as the former inviscid counterpart. A series of experimental studies under various conditions was made in the LENS II wind tunnel, and the RUN28 [26] case is chosen with inflow conditions of $Ma_\infty = 9.59, Re = 1.44 \times 10^5/m, T_\infty = 185.6$ K, and $T_W = 185.6$ K. The temporal algorithm adopts the LU-SGS scheme.

The geometries involved in the above cases are illustrated in subsequent result figures in Section 4.2 and thereby will not be repeated.

## 4.2 Numerical Results

In this subsection, we report tests of WENO3-$Z_{ES4}$ for the cases in Section 4.1, as well as comparative studies carried out using results provided by WENO3-Z [4–6], -F3 [10], and -$Z_{ES3}$ [15].

### 4.2.1 1D Scalar Equation

To validate the numerical accuracy of WENO3-$Z_{ES4}$, computations are performed using the initial condition given by Eq. (22). The $L_\infty$-norm errors and corresponding convergence rates for each scheme are tabulated in Table 2. The outcomes show that WENO3-Z fails to attain third-order convergence, while WENO-F3 exhibits improved accuracy but still falls short of third-order. In contrast, both WENO3-$Z_{ES3}$ and the proposed WENO3-$Z_{ES4}$ achieve third-order convergence. As indicated in Section 2, the failure of WENO-F3 to attain third-order accuracy in the $L_\infty$-norm originates from the occurrence of $CP_1$ at a half-node ($\lambda = 1/2$) at $t = 2$. Our theoretical investigations [13–15] established conclusively that all WENO3-Z improvements not grounded on ACPL inherently fail to achieve third-order in the $L_\infty$-norm in this case.

Table 2: $L_\infty$-norm errors and orders of WENO3-Z, WENO-F3, WENO3-$Z_{ES3}$, and -$Z_{ES4}$ by using equation of 1D scalar advection with initial condition Eq. (22) at $t = 2$ and CFL=0.25.

| N | Δt | WENO3-Z | | WENO-F3 | |
|---|---|---|---|---|---|
| | | $L_\infty$ −error | $L_\infty$ −order | $L_\infty$ −error | $L_\infty$ −order |
| 10 | 0.05 | 3.4899E-01 | - | 2.5309E-01 | - |
| 20 | 0.025 | 1.5189E-01 | 1.200 | 5.2671E-02 | 2.265 |
| 40 | 0.0125 | 6.3551E-02 | 1.257 | 7.1306E-03 | 2.885 |
| 80 | 0.00625 | 2.6923E-02 | 1.239 | 1.0222E-03 | 2.802 |
| 160 | 0.003125 | 1.0776E-02 | 1.321 | 1.6505E-04 | 2.631 |
| 320 | 0.0015625 | 4.0220E-03 | 1.422 | 3.0106E-05 | 2.455 |
| 640 | 0.00078125 | 1.4706E-03 | 1.451 | 7.3682E-06 | 2.031 |
| N | Δt | WENO3-$Z_{ES3}$ | | WENO3-$Z_{ES4}$ | |
| | | $L_\infty$ −error | $L_\infty$ −order | $L_\infty$ −error | $L_\infty$ −order |
| 10 | 0.05 | 2.9232E-01 | - | 2.1708E-01 | - |
| 20 | 0.025 | 7.6272E-02 | 1.938 | 4.6008E-02 | 2.238 |
| 40 | 0.0125 | 1.5602E-02 | 2.289 | 7.5831E-03 | 2.601 |
| 80 | 0.00625 | 8.0086E-03 | 0.962 | 1.0388E-03 | 2.868 |
| 160 | 0.003125 | 1.2820E-04 | 5.965 | 1.2814E-04 | 3.019 |
| 320 | 0.0015625 | 1.6035E-05 | 2.999 | 1.6035E-05 | 2.998 |
| 640 | 0.00078125 | 2.0047E-06 | 3.000 | 2.0047E-06 | 3.000 |

4.2.2 1D Euler Equations

1) Shu–Osher problem

To assess the resolution differences among the schemes, we choose a grid with 240 points rather than 400–600 points. The density distributions shown in Fig. 2 indicate that WENO3-Z has the poorest resolution performance at both peaks and valleys. WENO3-F3 ranks third in resolution performance, WENO3-$Z_{ES3}$ ranks second, and WENO3-$Z_{ES4}$ achieves the best resolution, obviously outperforming the other schemes in resolving all peaks and valleys. Notably, WENO3-$Z_{ES4}$ successfully resolves the secondary peak (marked by the dashed box), demonstrating its superior resolution. To the best of our knowledge, such secondary peaks and valleys are barely resolved on 240 grid points by typical third-order schemes.

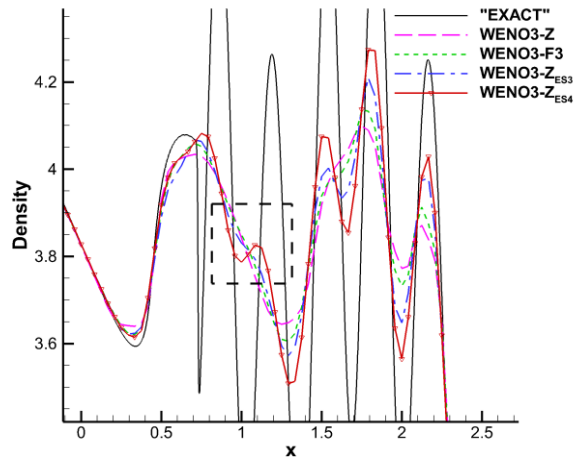

Figure 2: Density distributions for Shu–Osher problem at $t = 1.8$ on 240 grid points by using WENO3-Z, WENO-F3, WENO3-$Z_{ES3}$, and -$Z_{ES4}$.

2) Blast wave

The density distributions are shown in Fig. 3, demonstrating that all the schemes pass the test. Among them, WENO3-$Z_{ES4}$ has high resolution, particularly evident near the trough at $x = 0.746$. The resolution of WENO3-$Z_{ES3}$ is slightly lower than that of WENO3-$Z_{ES4}$ but marginally higher than those of WENO3-Z and WENO-F3.

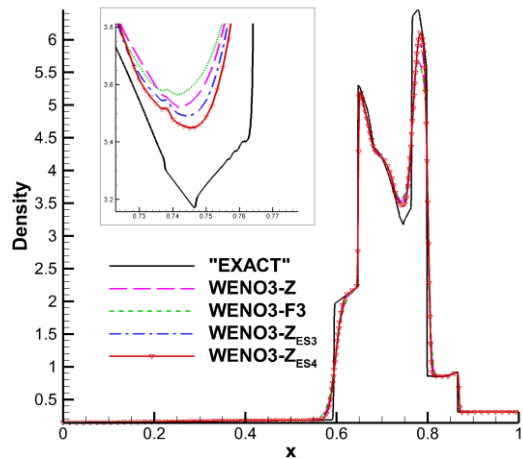

Figure 3: Density distributions for blast wave at $t = 0.038$ on 600 grid points by using WENO3-Z, WENO-F3, WENO3-$Z_{ES3}$, and -$Z_{ES4}$.

3) Strong shockwave

The density distributions are presented in Fig. 4. All the schemes complete the computations and confirm their numerical robustness for solving for a strong shockwave with a high pressure ratio. Local magnification of the results reveals that WENO3-$Z_{ES4}$ exhibits closer proximity to the theoretical solution at the top platform compared to the other schemes.

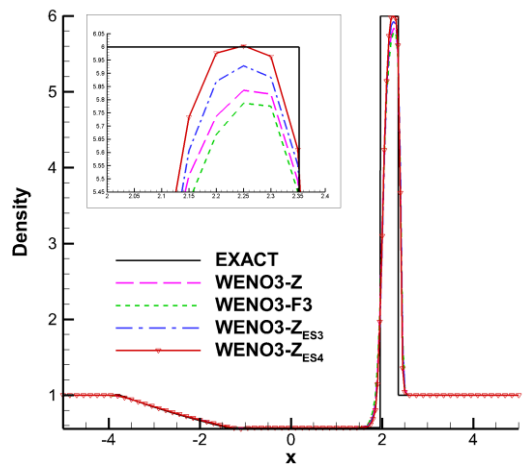

Figure 4: Density distributions for strong shockwave at $t = 0.01$ on 200 grid points with initial ratio of PR = $10^6$ by using WENO3-Z, WENO-F3, WENO3-$Z_{ES3}$, and -$Z_{ES4}$.

4.2.3  2D Euler Equations

1) 2D Riemann problem

Fig. 5 shows a comparison of the density contours for this problem. Although WENO3-Z produces smooth profiles, the absence of instabilities along the slip line implies its insufficient resolution. WENO-F3 demonstrates markable resolution improvements over WENO3-Z, particularly in resolving finer structures in the domain marked by the blue rectangle. However, in

Fig. 5(a), the pronounced asymmetry of structures therein implies compromised robustness. In contrast, WENO3-$Z_{ES3}$ resolves regular roll-ups along the slip line without significant numerical noise and exhibits enhanced structural symmetry within the domain, which confirms its balanced resolution and robustness. Compared to WENO3-$Z_{ES3}$, WENO3-$Z_{ES4}$ achieves analogous resolution but with more structural symmetry in the blue region. Furthermore, the density contours in the magnified window corresponding to those in the red box reveal that WENO3-$Z_{ES4}$ resolves not only the primary shear layer roll-up but also secondary structures, underscoring its superior resolution relative to the other schemes.

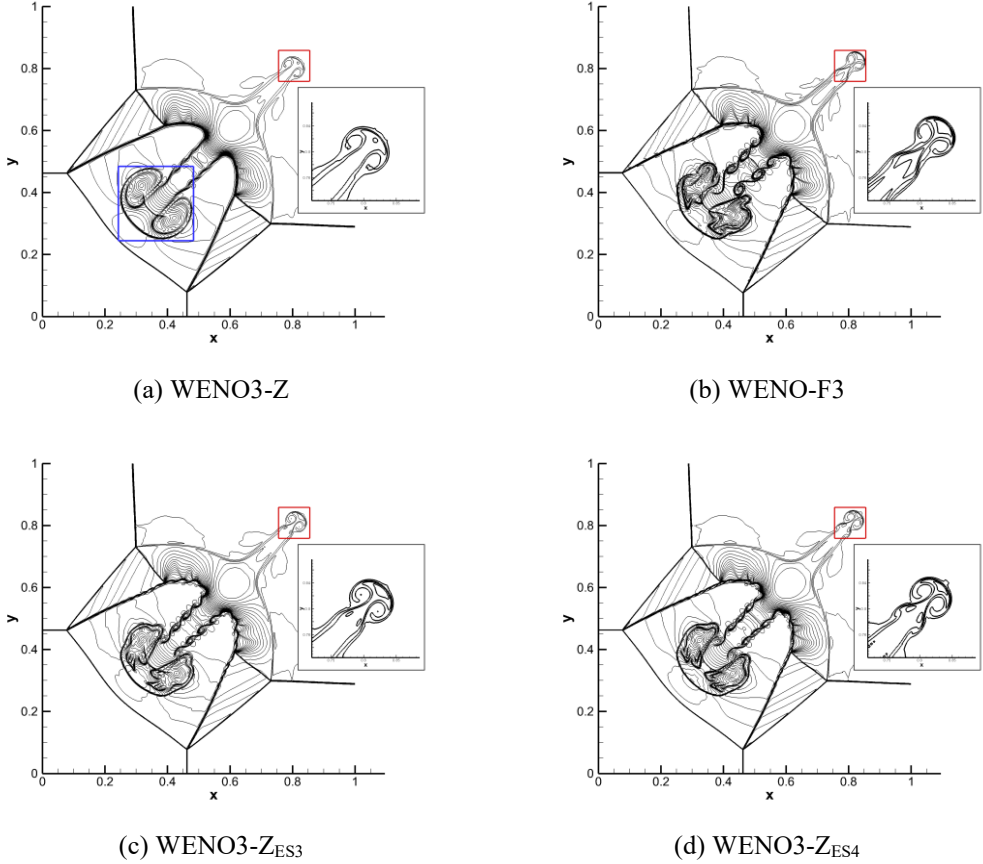

(a) WENO3-Z  (b) WENO-F3

(c) WENO3-$Z_{ES3}$  (d) WENO3-$Z_{ES4}$

Figure 5: Density contours of 2D Riemann problem by using WENO3-Z, WENO-F3, WENO3-$Z_{ES3}$, and -$Z_{ES4}$ on a $960 \times 960$ grid at $t = 0.8$ with $\Delta t = 0.0001$ (40 contours from 0.14 to 1.7).

2) Double Mach reflection

The density contours presented in Fig. 6 show the accomplishments of the four computational schemes. While unstable structures along the density slip line remain indistinct by WENO3-Z and WENO-F3, both WENO3-$Z_{ES3}$ and -$Z_{ES4}$ resolve these features. Notably, the latter gives a keen description of the structure in the blue rectangle in Fig. 6(d). For quantitative analysis, density distributions along $y = 0.06$ [marked by the dashed line in Fig. 6(a)] are provided in Fig. 7. All the schemes exhibit nearly the same discontinuities near $x \approx 2.78$, which indicates their consistency in solving for the primary shock (see Fig. 6). Fig. 7 further reveals density troughs and peaks around $x \approx 2.6 - 2.7$, corresponding to the central structure and trailing shock (indicated in the blue box in Fig. 6). Quantitative comparisons demonstrate that WENO3-$Z_{ES4}$ yields the lowest density of the central structure, followed sequentially by WENO3-$Z_{ES3}$, WENO-F3, and WENO3-

Z. Similarly, the magnitudes of peak density descend in the order of WENO3-$Z_{ES4}$, -$Z_{ES3}$, WENO-F3, and WENO3-Z. The regularity shown by the amplitudes and locations of the peaks and troughs by the four schemes suggests that WENO3-$Z_{ES4}$ has the highest resolution and least numerical dissipation.

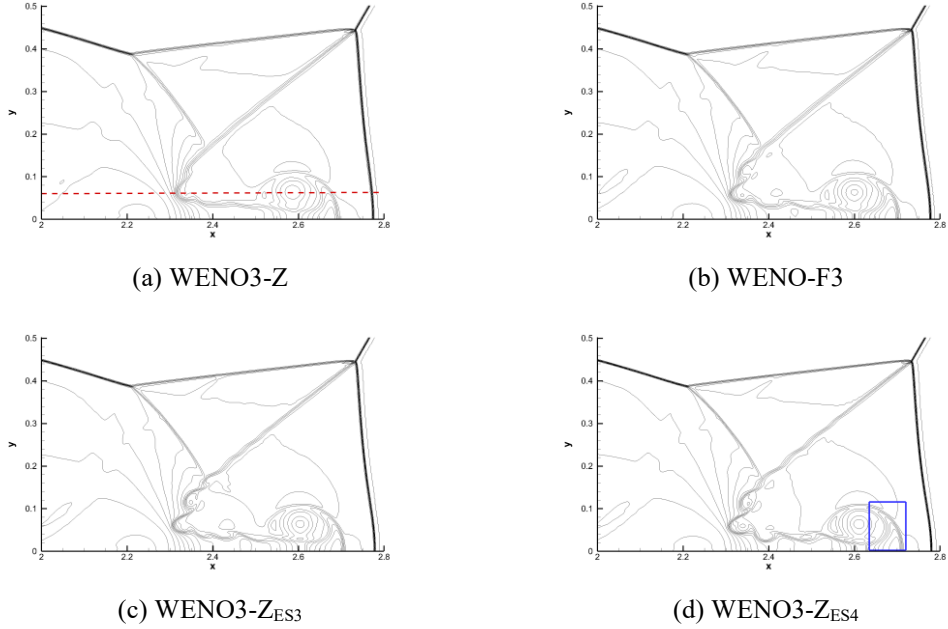

(a) WENO3-Z  (b) WENO-F3

(c) WENO3-$Z_{ES3}$  (d) WENO3-$Z_{ES4}$

Figure 6: Density contours of double Mach reflection by using WENO3-Z, WENO-F3, WENO3-$Z_{ES3}$, and -$Z_{ES4}$ on a $1920 \times 480$ grid at $t = 0.2$ with $\Delta t = 0.0001$ (33 contours from 1.4 to 24).

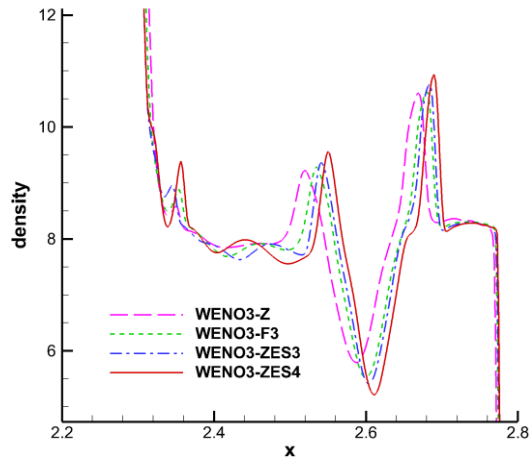

Figure 7: Density distributions on $y = 0.06$ of double Mach reflection by using WENO3-Z, WENO-F3, WENO3-$Z_{ES3}$, and -$Z_{ES4}$ on a $1920 \times 480$ grid at $t = 0.2$ and $\Delta t = 0.0001$.

3) Shock–bubble interaction

Density contours for this problem are shown comparatively in Fig. 8, with all simulations completed successfully. At $t = 10$, the gas bubble splits into two distinct regions because of shockwave interaction and initiates instability onset, allowing vortex roll-ups to serve as resolution metrics. Specifically, WENO3-Z produces the smoothest yet least-resolved results, while WENO-F3 demonstrates marginally improved resolution. Both WENO3-$Z_{ES3}$ and -$Z_{ES4}$ exhibit comparable vortex roll-ups, outperforming WENO-F3 in resolution. Comparatively, WENO3-$Z_{ES3}$ resolves the vortex structures of larger scale whereas WENO3-$Z_{ES4}$ yields the finer instabilities at the bubble connection.

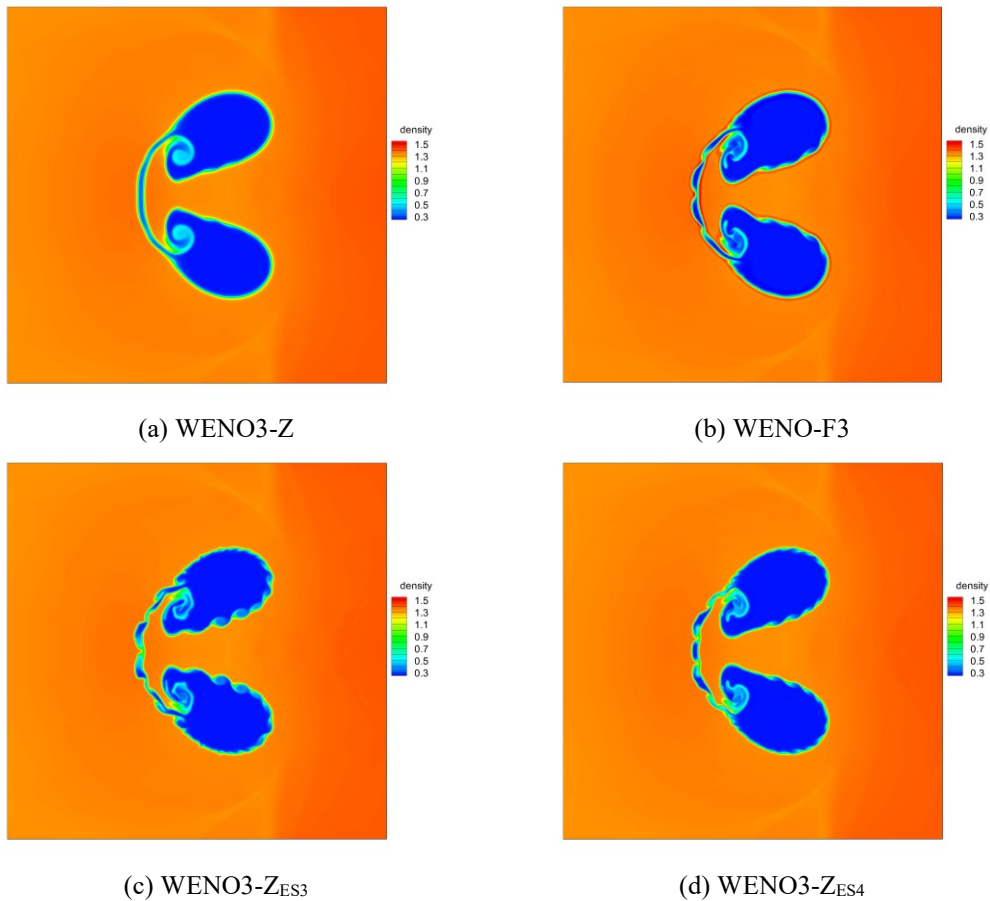

(a) WENO3-Z  (b) WENO-F3

(c) WENO3-$Z_{ES3}$  (d) WENO3-$Z_{ES4}$

Figure 8: Cloud maps of density contours from WENO3-Z, WENO-F3, WENO3-$Z_{ES3}$, and -$Z_{ES4}$ for shock–bubble interaction on a $1000 \times 1000$ grid at $t = 10$ with $\Delta t = 0.001$.

4) Mach-3 tunnel with a step

All schemes except WENO3-$Z_{ES3}$ completed the computation as shown in Fig. 9, implying the superior robustness of WENO3-$Z_{ES4}$ over WENO3-$Z_{ES3}$ in this problem. Despite the fact that the slip lines are all resolved downstream of the triple point, differences are observed in representing their instabilities. WENO-F3 and WENO3-$Z_{ES4}$ exhibit higher resolution than WENO3-Z. Furthermore, the WENO3-$Z_{ES4}$ results display more-pronounced vortex roll-up ahead of the reflected shock and a smoother post-leftward-shock flow field.

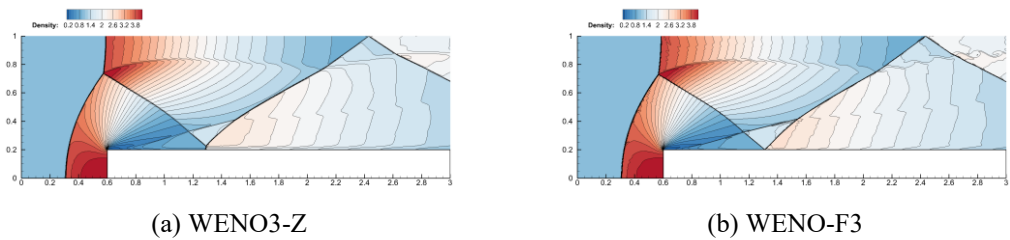

(a) WENO-Z  (b) WENO-F3

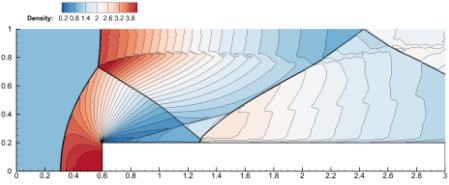

(c) WENO3-$Z_{ES4}$

Figure 9: Density contours from WENO3-Z, WENO-F3, and WENO3-$Z_{ES4}$ for a Mach-3 tunnel with a step on a $1200 \times 400$ grid at $t = 6$ and $\Delta t = 0.001$ (27 contours from 0.2 to 4.1).

5) Inviscid sharp-double-cone flow at $Ma_\infty = 9.59$

Mach-number contours given by the different schemes are shown comparatively in Fig. 10. The results of WENO3-Z exhibit minor oscillations in the upstream of the bow shock, while WENO-F3 gives the most pronounced fluctuations. In contrast, both WENO3-$Z_{ES3}$ and -$Z_{ES4}$ yield smoother shock profiles with less oscillations. A detailed examination of Fig. 10(b) reveals that the Mach-number oscillations in the perturbed region exceed the freestream value of 9.59. Therefore, the suppressed oscillations observed with WENO3-$Z_{ES3}$ and -$Z_{ES4}$ indicate their numerical robustness.

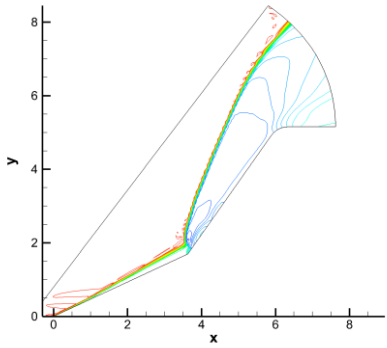  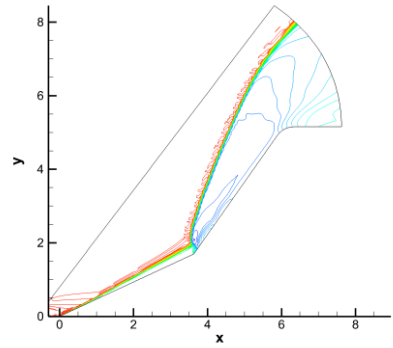

(a) WENO3-Z  (b) WENO-F3

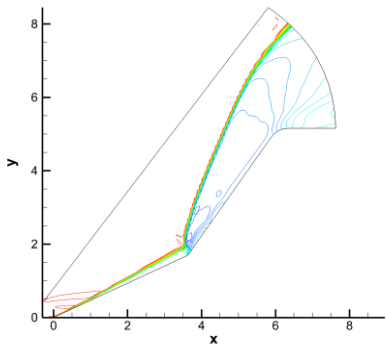  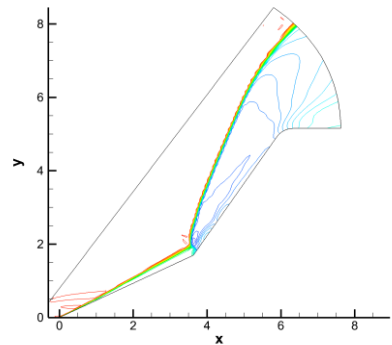

(c) WENO3-$Z_{ES3}$  (d) WENO3-$Z_{ES4}$

Figure 10: Mach-number contours of inviscid sharp-double-cone flow at $Ma = 9.59$ by using WENO3-Z, WENO-F3, WENO3-$Z_{ES3}$, and -$Z_{ES4}$ on a $204 \times 48$ gird at $t = 250$ with $\Delta t = 0.005$ (38 contours from 0 to 9.591).

6) Mach-2000 jet-flow problem

All the schemes successfully completed the computation, and their density contours in logarithmic scale are shown comparatively in Fig. 11. The overall flow structures align closely with those by Zhang and Shu [24], validating the numerical robustness of all the schemes for this problem.

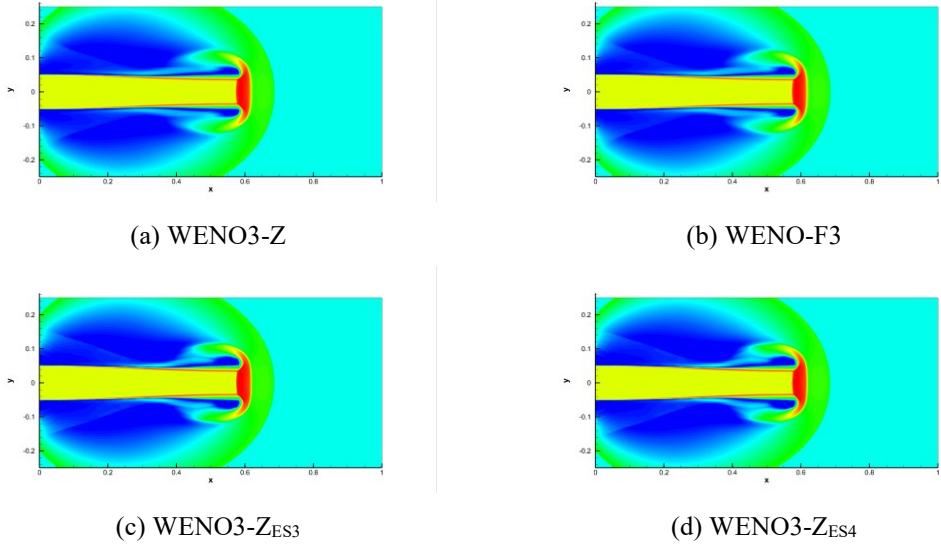

(a) WENO3-Z  (b) WENO-F3

(c) WENO3-$Z_{ES3}$  (d) WENO3-$Z_{ES4}$

Figure 11: Cloud maps of density contours in logarithmic-scale from WENO3-Z, WENO-F3, WENO3-$Z_{ES3}$, and -$Z_{ES4}$ for Mach-2000 jet flow on an $800 \times 400$ grid at $t = 10^{-3}$ with $\Delta t = 3 \times 10^{-7}$.

4.2.4 Navier–Stokes Equations

1) Flat-plate shock–boundary-layer interaction at $Ma_\infty = 2$

All four schemes completed the computation, and their pressure contours and streamline distributions are shown in Fig. 12. WENO3-Z yields the smallest separation bubble, while WENO3-$Z_{ES4}$ produces a slightly larger bubble than does WENO3-$Z_{ES3}$, comparable in size to that with WENO-F3. Different dissipation extents of the schemes are observed in the wave system after the bubble. For quantitative comparison, pressure distributions along $y = 0.405$ [marked by the dashed line in Fig. 12(a)] are shown in Fig. 13. As can be seen, WENO3-Z exhibits lower resolution at peaks and troughs compared to WENO3-$Z_{ES3}$, whereas WENO3-$Z_{ES4}$ outperforms WENO3-$Z_{ES3}$ but slightly lags behind WENO-F3. Overall, WENO3-$Z_{ES4}$ and WENO-F3 demonstrate better resolution and less dissipation.

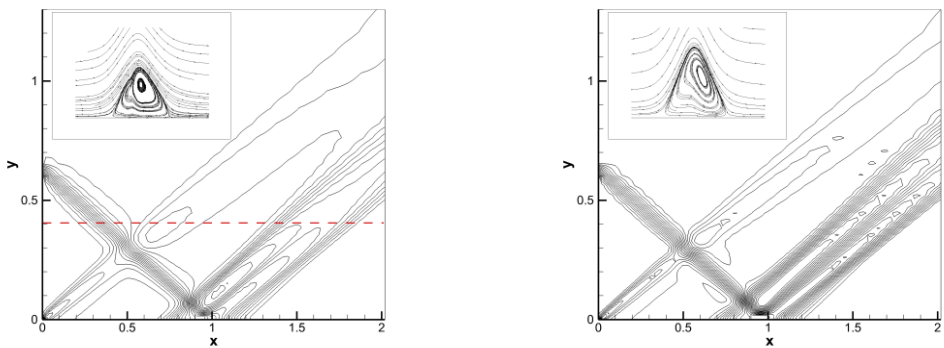

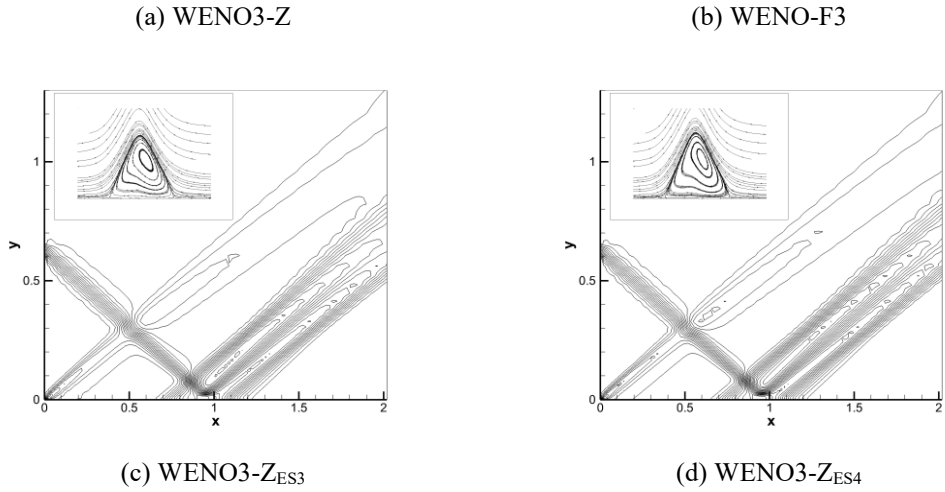

(a) WENO3-Z  (b) WENO-F3

(c) WENO3-$Z_{ES3}$  (d) WENO3-$Z_{ES4}$

Figure 12: Pressure contours and streamlines from WENO3-Z, WENO-F3, WENO3-$Z_{ES3}$, and -$Z_{ES4}$ for flat-plate shockwave–boundary-layer interaction at $Ma = 2$ on a $103 \times 122$ grid at $t = 40$ (28 contours from 0.18 to 0.245).

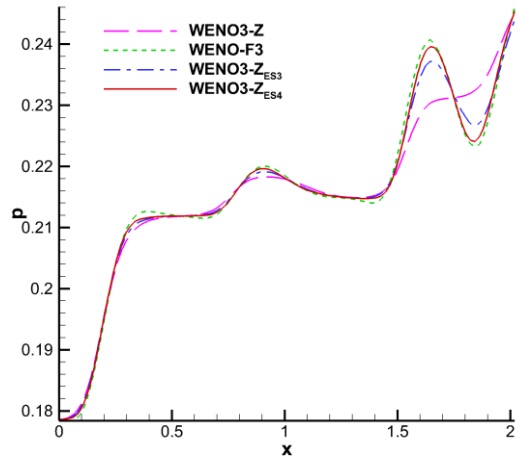

Figure 13: Pressure distributions on $y = 0.405$ for flat-plate shockwave–boundary-layer interaction at $Ma = 2$ by using WENO3-Z, WENO-F3, WENO3-$Z_{ES3}$, and -$Z_{ES4}$ on a $103 \times 122$ grid at $t = 40$.

2) Viscous shock-tube problem

Density contours obtained using the four schemes are displayed in Fig. 14. As is known, the triple-point location and adjacent shock structures are indicators of numerical resolution. Comparatively, WENO3-Z exhibits markedly lower resolution than do the other three schemes. To compare the vortex resolution ability, we analyze the one near $x \approx 0.825$. A blue line with $y = 0.14$ is drawn to indicate the vortex apex, which is chosen so that it aligns tangentially with the apex by WENO3-$Z_{ES4}$. The results show that WENO3-$Z_{ES4}$ achieves a vortex height comparable to that by WENO-F3 and greater than that by both WENO3-$Z_{ES3}$ and WENO3-Z, with WENO3-Z producing the lowest apex. Also, the vortex major axis (plotted in red dashed lines in the figure) is chosen as another resolution indicator by using its tilt angle. The measurements reveal the following sequence: WENO3-$Z_{ES4}$ gives the maximum angle (54.74°), followed by WENO3-$Z_{ES3}$ (44.18°), WENO-F3 (43.00°), and WENO3-Z (0.00°). The outcomes for vortex height and tilt angle of its major axis indicate conclusively that WENO3-$Z_{ES4}$ achieves an optimal resolution.

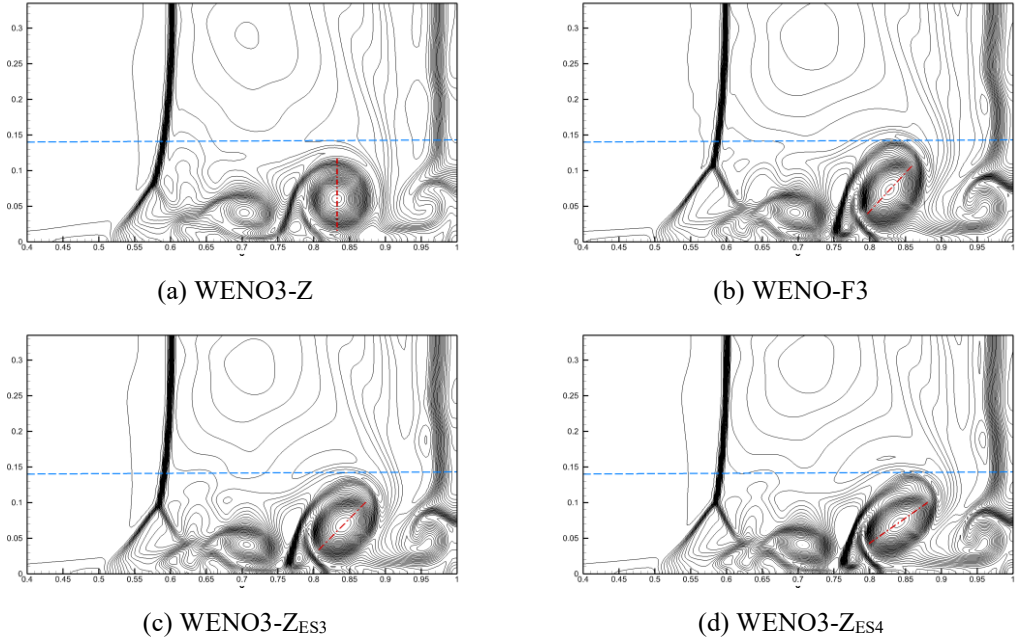

(a) WENO3-Z

(b) WENO-F3

(c) WENO3-$Z_{ES3}$

(d) WENO3-$Z_{ES4}$

Figure 14: Density contours from WENO3-Z, WENO-F3, WENO3-$Z_{ES3}$, and -$Z_{ES4}$ for viscous shock-tube problem on a $300 \times 300$ grid at $t = 1.0$ with $\Delta t = 0.004167$ (41 contours from 20 to 110).

3) Viscous sharp-double-cone flow at $Ma_\infty = 9.59$

Among the four schemes, WENO-F3 failed to complete the computation, indicating its insufficient robustness in this problem. For the three schemes that completed the computation, Fig. 15 shows their schlieren images using nondimensional density gradient ($|\nabla \rho / \nabla \rho_\infty|$) as well as wall heat flux distributions. To visualize the separation structures, streamlines within the separation zone are plotted. The results show that WENO3-Z produces the smallest separation vortex, while WENO3-$Z_{ES3}$ generates larger vortices than does WENO3-Z but slightly smaller than does WENO3-$Z_{ES4}$. Both WENO3-$Z_{ES4}$ and -$Z_{ES3}$ succeed in resolving secondary vortices. The wall heat flux in Fig. 15(d) reveals that WENO3-Z yields a separation obviously smaller than that of experiment [26], whereas the results of WENO3-$Z_{ES4}$ and -$Z_{ES3}$ exhibit reasonable agreement with the experimental data.

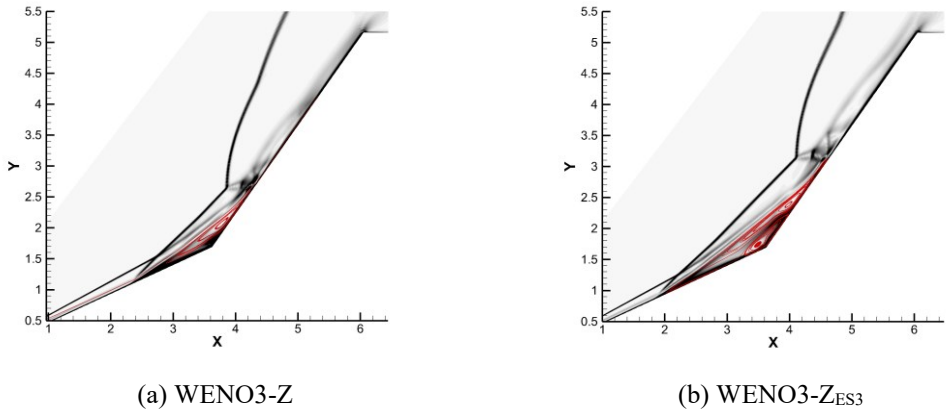

(a) WENO3-Z

(b) WENO3-$Z_{ES3}$

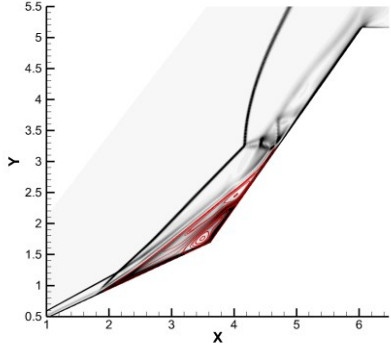 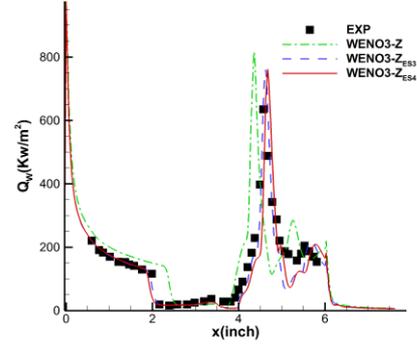

(c) WENO3-$Z_{ES4}$  (d) Distributions of wall heat- flux

Figure 15: (a)–(c) Contours of $|\nabla\rho/\nabla\rho_\infty|$ together with streamlines at separation corner of viscous sharp-double-cone flow by using WENO3-Z, -$Z_{ES3}$, and -$Z_{ES4}$; (d) distributions of wall heat flux by the same schemes with comparison against experimental results.

Finally, to summarize the robustness of the schemes, Table 3 gives their completion statuses for the various test cases, where "o" indicates accomplishment while "•" denotes failure. For simplicity, abbreviations are applied as follows: RM for the 2D Riemann problem, DMR for double Mach reflection, SBI for shock–bubble interaction, M3TS for the Mach-3 tunnel with a step, ISDC for inviscid sharp-double-cone flow, Jet$_{M2000}$ for the Mach-2000 jet-flow problem, SBLI for flat-plate shock–boundary-layer interaction, VST for the viscous shock tube, and SDC for viscous sharp-double-cone flow. The results show that the newly developed WENO3-$Z_{ES4}$ completed all the test cases, demonstrating its robustness while maintaining high resolution.

Table 3: Completion statuses of different schemes on different problems.

| Case<br>Scheme | RM | DMR | SBI | M3TS | ISDC | Jet$_{M2000}$ | SBLI | VST | SDC |
|---|---|---|---|---|---|---|---|---|---|
| WENO3-Z | o | o | o | o | o | o | o | o | o |
| WENO-F3 | o | o | o | o | o | o | o | o | • |
| WENO3-$Z_{ES3}$ | o | o | o | • | o | o | o | o | o |
| WENO3-$Z_{ES4}$ | o | o | o | o | o | o | o | o | o |

# 5 Conclusions

This study addressed the accuracy recovery at critical points of WENO3-Z and analyzed accuracy and resolution optimization. A method was developed to enhance the precision of nonlinear weights, via which a new WENO3-$Z_{ES4}$ was obtained that preserves the optimal order at $CP_1$ with $p=1$ in the nonlinear weights. Validating computations, comparisons, and analyses were presented, and the following conclusions are drawn.

1) An accuracy-optimization lemma was proposed, proving that the precision of normalized nonlinear weights can be improved under specified conditions ($\beta_k$ satisfy Eq. (12) and $C_{\alpha k}$ take the same value), through which the overall accuracy of the scheme can be elevated. Thus, an approach was established to enhance the accuracy order of schemes.

2) Analysis was conducted on the numerical resolution enhancement observed in computations when increasing $C_\beta$, leading to the proposal of Lemma 3.2, which shows that larger $C_\beta$ increases the relative importance of the discontinuous stencil, thereby causing the resolution enhancement.

3) Based on Lemma 3.1, $\beta_k^*$ and $\tau_4$ were constructed on extended stencils to satisfy the specified conditions. By integrating Lemmas 2.2 and 3.2 and via numerical experiments, the optimal values for $C_\beta$ and $C_\alpha$ were determined. The endeavor yielded the new scale-independent WENO3-$Z_{ES4}$ with optimal order recovered at $CP_1$.

4) Validating tests indicated that WENO3-$Z_{ES4}$ achieves a favorable balance between resolution and robustness. Specifically, in the Shu–Osher problem on 240 grid points, the scheme successfully resolved the density peak and valley of the second wave; in the 2D Riemann problem, it captured secondary structures in the head jets; in double Mach reflection and inviscid sharp-double-cone problems, the scheme delivered results with less numerical oscillation.

## Appendix

**Lemma 3.1:** For normalized nonlinear weights $\omega_k = \alpha_k / \sum_l \alpha_l$ where $\alpha_k = d_k \left(1 + C_{\alpha k} \left(\frac{\tau}{\beta_k}\right)^p\right)$ and $\tau = O(\Delta x^m)$, suppose that the accuracy relationship of $\beta_k$ satisfies $\beta_k = \sum_{l'=n_1}^{n_2-1} a_{l'} \Delta x^{l'} + b_k \Delta x^{n_2} + O(\Delta x^{n_2+1})$, where $b_k$ varies across different $k$, $n_2 > n_1 \geq 1$, and typically $m > n_1$. If $C_{\alpha k}$ take the same value for different $\alpha_k$, then $\omega_k = d_k \left(1 + O(\Delta x^{p(m-n_1)}) \times O(\Delta x^{n_2-n_1})\right)$.

**Proof:**

First, we show that $\left(\frac{\beta_k}{\beta_l}\right)^p = 1 + O(\Delta x^{n_2-n_1})$:

$$\frac{\beta_k}{\beta_l} = \frac{\sum_{l'=n_1}^{n_2-1} a_{l'} \Delta x^{l'} + b_k \Delta x^{n_2} + O(\Delta x^{n_2+1})}{\sum_{l'=n_1}^{n_2-1} a_{l'} \Delta x^{l'} + b_l \Delta x^{n_2} + O(\Delta x^{n_2+1})} = 1 + \frac{(b_k - b_l) \Delta x^{n_2} + O(\Delta x^{n_2+1})}{\sum_{l'=n_1}^{n_2-1} a_{l'} \Delta x^{l'} + b_l \Delta x^{n_2} + O(\Delta x^{n_2+1})}$$

$$= 1 + \Delta x^{n_2-n_1} \frac{(b_k - b_l)\Delta x^{n_1} + O(\Delta x^{n_1+1})}{\sum_{l'=n_1}^{n_2-1} a_{l'} \Delta x^{l'} + b_l \Delta x^{n_2} + O(\Delta x^{n_2+1})}$$

$$= 1 + \Delta x^{n_2-n_1} \frac{(b_k - b_l)\Delta x^{n_1}}{\sum_{l'=n_1}^{n_2-1} a_{l'} \Delta x^{l'} + b_l \Delta x^{n_2} + O(\Delta x^{n_2+1})}$$

$$= 1 + (b_k - b_l) \Delta x^{n_2-n_1} \frac{1 + O(\Delta x)}{\sum_{l'=0}^{n_2-n_1-1} a_{l'} \Delta x^{l'} + b_l \Delta x^{n_2-n_1} + O(\Delta x^{n_2+1})}$$

$$= 1 + (b_k - b_l) \Delta x^{n_2-n_1} (1 + O(\Delta x))$$

$$= 1 + O(\Delta x^{n_2-n_1}).$$

Thus, $\left(\frac{\beta_k}{\beta_l}\right)^p = 1 + O(\Delta x^{n_2-n_1})$. Next, we derive the accuracy relationship of $\omega_k$:

$$\omega_k = \frac{\alpha_k}{\sum_l \alpha_l} = \frac{d_k \left(1 + C_{\alpha k}\left(\frac{\tau}{\beta_k}\right)^p\right)}{1 + \sum_{l=0}^{r} d_l C_{\alpha l}\left(\frac{\tau}{\beta_l}\right)^p} = d_k \frac{\left(\frac{\beta_k}{\tau}\right)^p + C_{\alpha k}\left(\frac{\tau}{\tau \beta_k}\beta_k\right)^p}{\left(\frac{\beta_k}{\tau}\right)^p + \sum_{l=0}^{r} d_l C_{\alpha l}\left(\frac{\beta_k}{\beta_l}\right)^p}.$$

Substituting $\left(\frac{\beta_k}{\beta_l}\right)^p = 1 + O(\Delta x^{n_2-n_1})$, we have

$$\omega_k = d_k \frac{C_{\alpha k} + \left(\frac{\beta_k}{\tau}\right)^p}{\left(\frac{\beta_k}{\tau}\right)^p + \sum_{l=0}^{r} d_l C_{\alpha l}\left(\frac{\beta_k}{\beta_l}\right)^p} = d_k \frac{C_{\alpha k} + \left(\frac{\beta_k}{\tau}\right)^p}{\left(\frac{\beta_k}{\tau}\right)^p + \sum_{l=0}^{r} d_l C_{\alpha l}\left(1 + O(\Delta x^{n_2-n_1})\right)}.$$

Applying $\frac{a}{a+x} = 1 - \frac{1}{a}x + O(x^2)$, we have

$$\omega_k = d_k \frac{C_{\alpha k}+(\frac{\beta_k}{\tau})^p}{(\frac{\beta_k}{\tau})^p+C_{\alpha k}-C_{\alpha k}+\sum_{l=0}^r d_l C_{\alpha l}+O(\Delta x^{n_2-n_1})}$$

$$= d_k\left(1-\frac{-C_{\alpha k}+\sum_{l=0}^r d_l C_{\alpha l}+O(\Delta x^{n_2-n_1})}{C_{\alpha k}+(\frac{\beta_k}{\tau})^p}\right).$$

Let $A = -C_{\alpha k} + \sum_{l=0}^r d_l C_{\alpha l} + O(\Delta x^{n_2-n_1})$. Using $\frac{a}{a+x} = 1 - \frac{1}{a}x + O(x^2)$ again and substituting $\tau = O(\Delta x^m)$, we have

$$\omega_k = d_k\left(1 - \frac{1}{C_{\alpha k}+(\frac{\beta_k}{\tau})^p}A\right) = d_k\left(1 - C_{\alpha k}(\frac{\tau}{\beta_k})^p \times \frac{1}{C_{\alpha k}(\frac{\tau}{\beta_k})^p+1} \times A/C_{\alpha k}\right)$$

$$= d_k\left(1 - C_{\alpha k}\left(\frac{\tau}{\beta_k}\right)^p \times \left(1 - C_{\alpha k}\left(\frac{\tau}{\beta_k}\right)^p\right) \times A/C_{\alpha k}\right)$$

$$= d_k\left(1 - A(\frac{\tau}{\beta_k})^p\left(1 - C_{\alpha k}(\frac{\tau}{\beta_k})^p\right)\right)$$

$$= d_k(1 + A \times O(\Delta x^{p(m-n_1)})).$$

CASE 1: If $-C_{\alpha k} + \sum_{l=0}^r d_l C_{\alpha l} \neq 0$, we have

$$\omega_k = d_k\left(1 - (-C_{\alpha k} + \sum_{l=0}^r d_l C_{\alpha l})O(\Delta x^{p(m-n_1)})\right) = d_k\left(1 + O(\Delta x^{p(m-n_1)})\right),$$

obtaining the conventional error term.

CASE 2: If $-C_{\alpha k} + \sum_{l=0}^r d_l C_{\alpha l} = 0$ (i.e., $C_{\alpha k} = C_{\alpha l}$), we have

$$A = O(\Delta x^{n_2-n_1}) \rightarrow \omega_k = d_k(1 + O(\Delta x^{n_2-n_1})O(\Delta x^{p(m-n_1)})).$$

The accuracy of the error term is improved by $(n_2 - n_1)$ orders compared with the conventional $O(\Delta x^{p(m-n_1)})$.

**Lemma 3.2:** Consider two candidates $S_C$ and $S_D$ corresponding to the same stencil, where "C" and "D" denote that the variable is smoother on $S_C$ than on $S_D$ or $\beta_{k,D} > \beta_{k,C}$. Suppose $\beta_k = (\delta_j^{(1)m})^2 + C_{\beta,i}(\delta_j^{(2)n})^2$, where $\delta_j^{(n)m}$ represent the $n$-th order derivative with $m$-th order accuracy at point $x_j$, and the subscript $i$ indicates the different evaluations of $C_{\beta,i}$. Likewise, $\delta_{j,D}^{(1)m} > \delta_{j,C}^{(1)m}$ and $\delta_{j,D}^{(2)n} > \delta_{j,C}^{(2)n}$. For $\alpha_k = d_k(1 + \tau/\beta_k)$ with $\tau > 0$ and $\delta_{j,D}^{(1)m^2}/\delta_{j,C}^{(1)m^2} > \delta_{j,D}^{(2)n^2}/\delta_{j,C}^{(2)n^2}$, then $[\omega_{k,D}/\omega_{k,C}]_{C_{\beta,1}} > [\omega_{k,D}/\omega_{k,C}]_{C_{\beta,2}}$ providing that $C_{\beta,1} > C_{\beta,2} > 0$.

**Proof**:

Observing that $[\frac{\omega_{k,D}}{\omega_{k,C}}]_{C_{\beta,i}} > [\frac{\alpha_{k,D}}{\alpha_{k,C}}]_{C_{\beta,i}} = \left(1 + \frac{\tau}{(\delta_{jD}^{(1)m})^2+x(\delta_{jD}^{(2)n})^2}\right)\Big/\left(1 + \frac{\tau}{(\delta_{jC}^{(1)m})^2+x(\delta_{jC}^{(2)n})^2}\right)$, we consider the function $y(x) = \left(1 + \frac{\tau}{(\delta_{jD}^{(1)m})^2+x(\delta_{jD}^{(2)n})^2}\right)\Big/\left(1 + \frac{\tau}{(\delta_{jC}^{(1)m})^2+x(\delta_{jC}^{(2)n})^2}\right)$ and compute its first derivative:

$$y'(x) = \left(\left(-\left(\delta_{jC}^{(1)m}\right)^2 \left(\delta_{jD}^{(2)n}\right)^2 + \left(\delta_{jD}^{(1)m}\right)^2 \left(\delta_{jC}^{(2)n}\right)^2\right)\tau\right.$$
$$+ \left(-\left(\delta_{jC}^{(1)m}\right)^4 \left(\delta_{jD}^{(2)n}\right)^2 + \left(\delta_{jD}^{(1)m}\right)^4 \left(\delta_{jC}^{(2)n}\right)^2\right)$$
$$+ \left(-\left(\delta_{jC}^{(2)n}\right)^4 \left(\delta_{jD}^{(2)n}\right)^2 + \left(\delta_{jC}^{(2)n}\right)^2 \left(\delta_{jD}^{(2)n}\right)^4\right)x^2$$
$$+ \left(-2\left(\delta_{jC}^{(1)m}\right)^2 \left(\delta_{jC}^{(2)n}\right)^2 \left(\delta_{jD}^{(2)n}\right)^2 + 2\left(\delta_{jD}^{(1)m}\right)^2 \left(\delta_{jC}^{(2)n}\right)^2 \left(\delta_{jD}^{(2)n}\right)^2\right)x\right)$$
$$\left./\left(\left(\delta_{jC}^{(1)m}\right)^2 + \tau + x\left(\delta_{jC}^{(2)n}\right)^2\right)^2 \left(\left(\delta_{jD}^{(1)m}\right)^2 + x\left(\delta_{jD}^{(2)n}\right)^2\right)^2\right).$$

The denominator is strictly positive, and the numerator constitutes a quadratic function in $x$. Given $\delta_{j,D}^{(1)m} > \delta_{j,C}^{(1)m}$ and $\delta_{j,D}^{(2)n} > \delta_{j,C}^{(2)n}$, the following inequalities hold:

$$-\left(\delta_{jC}^{(2)n}\right)^4 \left(\delta_{jD}^{(2)n}\right)^2 + \left(\delta_{jC}^{(2)n}\right)^2 \left(\delta_{jD}^{(2)n}\right)^4 > 0,$$

$$\left(-2\left(\delta_{jC}^{(1)m}\right)^2 \left(\delta_{jC}^{(2)n}\right)^2 \left(\delta_{jD}^{(2)n}\right)^2 + 2\left(\delta_{jD}^{(1)m}\right)^2 \left(\delta_{jC}^{(2)n}\right)^2 \left(\delta_{jD}^{(2)n}\right)^2\right) > 0.$$

Thus, the linear and quadratic coefficients of the numerator are strictly positive. Consequently, the quadratic function's symmetry axis lies on the negative $x$-axis. When the constant term of the numerator satisfies non-negativity, i.e.,

$$\left(-\left(\delta_{jC}^{(1)m}\right)^2 \left(\delta_{jD}^{(2)n}\right)^2 + \left(\delta_{jD}^{(1)m}\right)^2 \left(\delta_{jC}^{(2)n}\right)^2\right)\tau \geq 0,$$

$$\left(-\left(\delta_{jC}^{(1)m}\right)^4 \left(\delta_{jD}^{(2)n}\right)^2 + \left(\delta_{jD}^{(1)m}\right)^4 \left(\delta_{jC}^{(2)n}\right)^2\right) \geq 0,$$

the first derivative $y'(x)$ remains positive for all $x > 0$. The two conditions reduce to $\frac{\delta_{j,D}^{(1)m^2}}{\delta_{j,C}^{(1)m^2}} > \frac{\delta_{j,D}^{(2)n^2}}{\delta_{j,C}^{(2)n^2}}$ and $\left(\frac{\delta_{j,D}^{(1)m^2}}{\delta_{j,C}^{(1)m^2}}\right)^2 > \frac{\delta_{j,D}^{(2)n^2}}{\delta_{j,C}^{(2)n^2}}$. Given $\delta_{j,D}^{(1)m} > \delta_{j,C}^{(1)m}$, then $\left(\frac{\delta_{j,D}^{(1)m^2}}{\delta_{j,C}^{(1)m^2}}\right)^2 > \frac{\delta_{j,D}^{(1)m^2}}{\delta_{j,C}^{(1)m^2}} > \frac{\delta_{j,D}^{(2)n^2}}{\delta_{j,C}^{(2)n^2}}$. Thus, the combined conditions reduce to $\frac{\delta_{j,D}^{(1)m^2}}{\delta_{j,C}^{(1)m^2}} > \frac{\delta_{j,D}^{(2)n^2}}{\delta_{j,C}^{(2)n^2}}$. For $a > b > 0$ under this condition $y(a) > y(b)$ is ensured. Consequently, when $C_{\beta,1} > C_{\beta,2} > 0$, the normalized weight ratio satisfies $\left[\frac{\omega_{k,D}}{\omega_{k,C}}\right]_{C_{\beta_1}} > \left[\frac{\omega_{k,D}}{\omega_{k,C}}\right]_{C_{\beta_2}}$.